\documentclass[a4paper,11pt,twoside]{article}
\usepackage{pst-fill,pst-grad}
\usepackage{textcomp}
\usepackage[english]{babel}
\usepackage[utf8x]{inputenc}
\usepackage{graphicx}
\usepackage{amsmath}
\usepackage{multicol}
\usepackage{float}
\usepackage{fancyhdr}
\usepackage[matrix,arrow,curve]{xy}
\usepackage{pstricks} 
\usepackage{amsmath,amsfonts,verbatim,afterpage,theorem,euscript,mathrsfs,amssymb}
\usepackage{amsfonts}
\usepackage{amssymb}
\usepackage{array}
\usepackage{dsfont}
\usepackage{hyperref}
\usepackage{enumitem}

\numberwithin{equation}{section}
\def \vu{\vec{u}}
\def \vv{\vec{v}}

\def \vU{\vec{U}}
\def \vB{\vec{B}}
\def \vA{\vec{A}}

\def \vn{\vec{\nabla}}
\def \vb{\vec{b}}
\def \vh{\vec{h}}
\def \vk{\vec{k}}
\def \vl{\vec{l}}

\def \vbe{\vec{\beta}}
\def \vga{\vec{\gamma}}

\def \vf{\vec{f}}
\def \vg{\vec{g}}
\def \vF{\vec{F}}
\def \vG{\vec{G}}
\def \vX{\vec{X}}
\def \vY{\vec{Y}}
\def \vZ{\vec{Z}}
\def \ep{\varepsilon}

\def \vphi{\vec{\varphi}}

\def \vpphi{\vec{\phi}}

\def \R{\mathbb{R}^{3}}

\def \Qr{Q_{\rho}}
\def \Qro{Q_{\rho_0}}
\def \Qrone{Q_{\rho_1}}
\def \Br{B_{\rho}}
\def \Bro{B_{\rho_0}}

\def \pat{\partial_t}

\newcommand{\norm}[1]{\lVert#1\rVert}

\newtheorem{Definition}{Definition}[section]
\newtheorem{Proposition}{Proposition}[section]
\newtheorem{Lemme}{Lemma}[section]
\newtheorem{Theoreme}{Theorem}

\newtheorem{Remarque}{Remark}[section]

\title{\bf  Regularity theory for the dissipative\\ solutions of the MHD equations}
\usepackage{authblk}
\author{Diego Chamorro\footnote{\emph{diego.chamorro@univ-evry.fr}} }
\author{Jiao He\footnote{\emph{jiao.he@univ-evry.fr}}}
\affil{\footnotesize LaMME, Univ. Evry, CNRS, Universit\'e Paris-Saclay, 91025, Evry, France.}

\begin{document}
\maketitle
\begin{scriptsize}
\abstract{We study here a new generalization of Caffarelli, Kohn and Nirenberg's partial regularity theory for weak solutions of the MHD equations. 
Indeed, in this framework some hypotheses on the pressure $P$ are usually asked (for example $P\in L^q_tL^1_x$ with $q>1$) and then local H\"older regularity,  in time and space variables, for weak solutions can be obtained over small neighborhoods.
By introducing the notion of dissipative solutions, we weaken the hypothesis on the pressure (we will only assume that $P\in \mathcal{D}'$) and we will obtain H\"older regularity in the space variable for weak solutions.}\\	

\noindent\textbf{Keywords:} MHD equations; dissipative solutions; partial regularity theory; Morrey spaces.\\
\textbf{MSC2020:} 76W05; 76D03; 35Q35.
\end{scriptsize}

\section{Introduction}
The main purpose of this work is to study the partial regularity theory for the incompressible 3D MagnetoHydroDynamic (MHD) equations which are given by the following system:
\begin{equation}\label{EquationMHDoriginal}
\begin{cases}
\partial_{t}\vU=\Delta \vU -(\vU\cdot\vn)\vU+(\vB\cdot\vn)\vB-\vn P+\vF,\quad div(\vU) = div(\vF)=0,\\[3mm]
\partial_{t}\vB=\Delta \vB -(\vU\cdot\vn)\vB+(\vB\cdot\vn)\vU+\vG,\quad div(\vB)=div(\vG)=0,\\[3mm]
\vU(0,x)=\vU_{0}(x), \; div(\vU_0)=0 \mbox{ and } \vB(0,x)=\vB_{0}(x),\; div(\vB_0)=0, \qquad x\in \mathbb{R}^3,
\end{cases}
\end{equation}
where $\vU, \vB:[0,T]\times \R\longrightarrow \R$ are two divergence-free vector fields which represent the velocity and the magnetic field, respectively, and the scalar function $P:[0,T]\times \R\longrightarrow \mathbb{R}$ stands for the pressure. The initial data $\vU_0, \vB_0: \R\longrightarrow \R$ and the external forces $\vF, \vG:[0,T]\times \R\longrightarrow \R$ are given. We will say that $(\vU, P, \vB)$ is a solution of the MHD equations if it solves the system (\ref{EquationMHDoriginal}) in some particular sense to be precised later on.\\

In the recent articles \cite{ChCHJ1} and \cite{ChCHJ2} we have studied two different regularity theories for the system (\ref{EquationMHDoriginal}) using parabolic Morrey spaces as a common framework. The aim of the current work is to mix in a very specific manner these two previous results to obtain a new class of solutions of (\ref{EquationMHDoriginal}) for which we can deduce more general regularity results.\\

Indeed, in the first article \cite{ChCHJ1} we obtained a \emph{local} regularity result for the MHD equations following some ideas of O'Leary \cite{OLeary} (which were originally stated for the classical 3D Navier-Stokes equations): if we assume that $\mathds{1}_{\Omega}\vU$ and $\mathds{1}_{\Omega}\vB$ belong to some suitable parabolic Morrey spaces $\mathcal{M}^{p,q}_{t,x}$ (see Section \ref{Section_Preliminaries} below for a definition of these functional spaces), where $\Omega$ is a bounded subset of $]0, +\infty[\times \mathbb{R}^3$ of the form
\begin{equation}\label{DefConjuntoOmega}
\Omega=]a,b[\times B(x_{0},r), \quad \mbox{with} \quad 0<a<b<+\infty, x_{0}\in \R \mbox{ and } 0<r<+\infty. 
\end{equation}
Then it is possible to obtain a gain of regularity in the space variable even though the pressure is a general object (\emph{i.e.} $P\in \mathcal{D}'$). The  result obtained in this framework is the following one:
\begin{Theoreme}[Local Regularity, \cite{ChCHJ1}]\label{Teo_SerrinMHD}
Let $\vU_{0}, \vB_{0}:\R\longrightarrow \R$ such that $\vU_{0}, \vB_{0}\in L^{2}(\R)$ and $div(\vU_{0})=div(\vB_{0})=0$ be two initial data and consider two external forces $\vF, \vG:[0,+\infty[\times\R\longrightarrow \R$ such that $\vF, \vG\in L^{2}([0,+\infty[, \dot{H}^{1}(\R))$. Assume that $\Omega$ is a bouned set of the form (\ref{DefConjuntoOmega}), that $P\in \mathcal{D}'(\Omega)$ and that $\vU, \vB:[0,+\infty[\times\R\longrightarrow \R$ are two vector fields that belong to the space
\begin{equation}\label{UBSolutionFaiblesMHD}
L^{\infty}(]a,b[, L^{2}(B(x_{0},r)))\cap L^{2}(]a,b[, \dot{H}^{1}(B(x_{0},r))),
\end{equation}
such that they satisfy the MHD equations (\ref{EquationMHDoriginal}) over the set $\Omega$.\\
		
\noindent If moreover we have the following local hypotheses
\begin{equation}\label{LocalHypo1}
\begin{cases}
\mathds{1}_{\Omega}\vU\in\mathcal{M}_{t,x}^{p_{0},q_{0}}(\mathbb{R}{\tiny \times}\R) & \mbox{ with } 2<p_{0}\leq q_{0}, 5<q_{0}<+\infty\\[3mm]
\mathds{1}_{\Omega}\vB\in\mathcal{M}_{t,x}^{p_{1},q_{1}}(\mathbb{R}\times\R) &\mbox{ with } 2<p_{1}\leq q_{1}, 5<q_{1}<+\infty, 
 \end{cases}
\end{equation} 
and $p_1 \leq p_0$, $q_1 \leq q_0$, then, for all $\alpha, \beta$ such that $a<\alpha<\beta<b$ and for all $\rho$ such that $0<\rho<r$, we have 
\begin{equation}\label{ConclusionSerrinMHD}
\vU\in L^{q_0}(]\alpha,\beta[, L^{q_0}(B(x_{0},\rho)))\quad \mbox{and}\quad \vB\in L^{q_1}(]\alpha,\beta[, L^{q_1}(B(x_{0},\rho))).
\end{equation}
\end{Theoreme}
Following  Serrin's theory \cite{Serrin1} (see also Section 13.2 of \cite{PGLR1}) it can be shown that the conclusion of this theorem implies that 
\begin{equation}\label{ConclusionSerrinMHD1}
\vU\in L^{\infty}(]\alpha',\beta'[, L^{\infty}(B(x_{0},\rho')))\quad \mbox{and}\quad \vB\in L^{\infty}(]\alpha',\beta'[, L^{\infty}(B(x_{0},\rho'))),
\end{equation}
for some $\alpha\leq \alpha'<\beta'\leq \beta$ and $0<\rho'\leq \rho$. 
\begin{Remarque}\label{Rem_RegulariteSerrin}
The points $(t,x)\in [0,+\infty[\times \mathbb{R}^3$ such that the previous condition (\ref{ConclusionSerrinMHD1}) holds for $\vU(t,x)$ and $\vB(t,x)$ are called \emph{regular points} as we can prove that in this case the regularity (in the space variable) of the solution $\vU$ and $\vB$ is actually driven by the regularity of the external forces $\vF$ and $\vG$: in the framework of Theorem \ref{Teo_SerrinMHD} we obtain that $\vU, \vB\in L^\infty_tH^2_x\cap L^2_t H^3_x$ over a subset of $\Omega$. See Theorem 13.1 p. 397 of \cite{PGLR1} for a proof of this fact in the setting of the classical Navier-Stokes equations which can be easily extended to the MHD equations.
\end{Remarque} 
These ideas are of course reminiscent from the work of Serrin \cite{Serrin1} (stated for the classical Navier-Stokes equations) which first considered the condition $\vU\in (L^p_tL^q_x)_{loc}$ with $\tfrac{2}{p}+\tfrac{3}{q}<1$ (the case $\tfrac{2}{p}+\tfrac{3}{q}=1$ was studied by  Struwe \cite{Struwe} and Takahashi \cite{Takahashi}). One important point of all these results is related to Serrin's counter-example: as we do not impose any particular information on the pressure $P$, we may lose regularity in the time variable since the temporal regularity is closely linked to the pressure (see Section 13.1 of \cite{PGLR1} for this particular point) and thus only space regularity can be obtained. Note also that the Morrey space assumption (\ref{LocalHypo1}) or the $(L^p_tL^q_x)_{loc}$ control with $\tfrac{2}{p}+\tfrac{3}{q}\leq 1$ are quite strong hypotheses that can not be deduced from a general weak Leray-type solution of (\ref{EquationMHDoriginal}).

Let us now mention that in the study of the local regularity theory for the MHD equations, Wu \cite{Wu1} generalized Ladyzhenskaya-Prodi-Serrin-type criteria by assuming conditions on both the velocity $\vU$ and the magnetic field $\vB$. He and Xin \cite{He1} gave a Serrin-type regularity condition which only including the integrability condition on the velocity $\vU$. More precisely, they shown that weak solutions $(\vU, \vB)$ are smooth if $\vU \in (L^p_tL^q_x)_{loc}$ with $\tfrac{2}{p}+\tfrac{3}{q}\leq 1$ and $q >3$, which reveals that the velocity field plays a more dominant role than the magnetic field does on the regularity of solutions to the MHD equations.  
Chen, Miao and Zhang \cite{ChenMZ} then proved that $(\vU, \vB)$ is regular under a refined Serrin-type regularity criterion in the framework of Besov spaces with negative index in terms of the velocity only. Later on, Kang and Lee \cite{KangLee} present a new interior regularity criteria for suitable weak solutions to the MHD equations by supposing that some scaled norm of the velocity is small and some scaled norm of the magnetic field is bounded. Their result was improved by Wang and Zhang in \cite{WangZhang} by removing the bounded assumptions on the magnetic field. \\

In the second article \cite{ChCHJ2} we studied the so-called \emph{partial} regularity theory for the MHD equations (\ref{EquationMHDoriginal}) using the same parabolic Morrey spaces setting as above. This setting was initially applied for the Navier-Stokes equations by Kukavica in \cite{Kukavica} and \cite{Kukavica08} in order to generalize the Caffarelli, Kohn and Nirenberg regularity criterion \cite{CKN}. The main idea of this theory relies in the use of the of ``suitable weak solutions" (first developed by Scheffer for the Navier-Stokes equations in \cite{Scheffer} and \cite{Scheffer1})  which is a class of weak Leray-type solution $(\vU, P, \vB)$ such that the distribution $\mu$ given by the expression
\begin{eqnarray}
\mu&=&-\partial_t(|\vU|^2 + |\vB|^2 )+ \Delta (|\vU|^2 + |\vB|^2 )  - 2 (|\vn \otimes \vU|^2 + |\vn \otimes \vB|^2 ) - 2 div (P \vU)   \label{Formula_MesureMHD}\\  
& &- div \left( (|\vU|^2 + |\vB|^2) \vU \right) + 2 div ((\vU \cdot \vB) \vB ) + 2 (\vF \cdot \vU + \vG \cdot\vB),\notag
\end{eqnarray}
is a non-negative locally finite measure on $\Omega$, where $\Omega\subset ]0,+\infty[\times \mathbb{R}^3$ is a bounded domain. Of course, in the formula above some conditions on the pressure $P$ must be imposed (usually $P\in L^q_tL^1_x$ with $q>1$) in order to make the quantity $div(P \vU)$ meaningful. Assuming moreover a particular behavior of the quantities $\vn\otimes \vU$ and $\vn\otimes \vB$ over a small neighborhood of points, then it is possible to obtain a gain of regularity on both variables, space and time. More precisely we have: 
\begin{Theoreme}[Partial Regularity, \cite{ChCHJ2}]\label{Teo_CKNMHD} 
Let $\Omega$ be a bounded domain of the form given in (\ref{DefConjuntoOmega}). Let $(\vU, P, \vB)$ be a weak solution on $\Omega$ of the MHD equations \eqref{EquationMHDoriginal}.
Assume that
\begin{itemize}
\item[1)] The vector $(\vU, \vB, P, \vF, \vG)$ satisfies the conditions 
$$\vU, \vB \in L^{\infty}_t L^{2}_x \cap L^{2}_t \dot{H}^{1}_x(\Omega),\qquad 
P \in L^{q_{0}}_{t,x}(\Omega) \mbox{ with } 1<q_0\leq \tfrac{3}{2},\qquad
\vF, \vG \in L^{\frac{10}{7}}_{t,x}(\Omega).$$
\item[2)] The solution $(\vU, P, \vB)$ is suitable in the sense that the distribution (\ref{Formula_MesureMHD}) is a non-negative locally finite measure on $\Omega$.
\item[3)] We have the following local information on $\vF$ and $\vG$: $\mathds{1}_{\Omega}\vF \in \mathcal{M}_{t,x}^{\frac{10}{7}, \tau_{a}}$ and $ \mathds{1}_{\Omega} \vG \in \mathcal{M}_{t,x}^{\frac{10}{7}, \tau_{b}}$ for some $\tau_{a}, \tau_b>\frac{5}{2-\alpha}$ with $0<\alpha<\frac{1}{3}$.
\end{itemize}
There exists a positive constant $\epsilon^{*}$ which depends only on $\tau_{a}$ and $\tau_b$ such that, if for some $\left(t_{0}, x_{0}\right) \in \Omega$, we have
\begin{equation}\label{HypothesePetitesseGrad}
\limsup _{r \rightarrow 0} \frac{1}{r} \iint_{]t_{0}-r^{2}, t_{0}+r^{2}[ \times B\left(x_{0}, r\right)}|\vn \otimes \vU|^{2} + |\vn \otimes \vB|^{2} dxds <\epsilon^{*},
\end{equation}
then $(\vU, \vB)$ is H\"older regular of exponent $\alpha$ (in the time variable and in the space variable) in a neighborhood of $\left(t_{0}, x_{0}\right)$ for some small $\alpha$ in the interval $0<\alpha<\frac{1}{3}$.
\end{Theoreme}

It is worth noting here that considerable efforts have been made to weaken as much as possible the hypotheses on the pressure. Indeed, in the original paper of Caffarelli, Kohn and Nirenberg \cite{CKN} which deals with the Navier-Stokes equations, the authors assumed that $P \in (L^{\frac{5}{4}}_{t,x})_{loc}$. Later on, based on a better estimate on the pressure due to Sohr and von Wahl \cite{Sohrvon}, Lin proposed a simplified proof in \cite{Lin} and assumed that $P \in (L^{\frac{3}{2}}_{t,x})_{loc}$. Afterwards, under a more natural condition on the pressure, Seregin and Sverak \cite{SereginSve} showed the regularity of suitable solutions to the Navier-Stokes equations. Vasseur \cite{Vasseur} gave a proof relying on a method introduced by De Giorgi for elliptic equations to show that $P \in (L^{q_0}_{t}L^{1}_{x})_{loc}$ with $q_0 >1$ is actually enough. 

Note also that in the setting of the MHD equations, He and Xin \cite{He2} studied uniform gradient estimates and extended the work of Caffarelli, Kohn and Nirenberg \cite{CKN} to the partial regularity theory of suitable weak solutions to the MHD equations under no condition for the magnetic field $\vB$. Recently, relying on the De Giorgi iteration developed by Vasseur for the Navier-Stokes eqautions, Jiu and Wang \cite{Jiu_Wang14} gave an alternative proof of the work of He and Xin \cite{He2} to show that one dimensional parabolic Hausdorff measure of
the possible singular points is zero.
In \cite{CaoWU}, Cao and Wu provided a regularity criteria suggesting that the directional derivative of the pressure $P$ is bounded in some Lebesgue spaces. We remark that all these results hold under some integrability conditions on the pressure.\\

In this article we are going one step further in the treatment of the pressure and following the ideas of \cite{CML} we will prove that a pressure $P\in \mathcal{D}'$ can be considered to obtain a more general version of Theorem \ref{Teo_CKNMHD}.  Indeed, although the \emph{local} and \emph{partial} regularity theories are quite different in spirit, it is possible to combine them to obtain some new results on the regularity of the solutions of the MHD equations (\ref{EquationMHDoriginal}). To be more precise, we will use the local regularity theory given in Theorem \ref{Teo_SerrinMHD} (which does not imposes conditions on the pressure $P$) in order to generalize the partial regularity theory stated in Theorem \ref{Teo_CKNMHD}: even though $P$ is a distribution, we can give a sense to the expression (\ref{Formula_MesureMHD}) above and with some additional mild assumptions we will be able to obtain a gain of regularity. \\

Remark that the use of the local regularity theory will not be direct and we need to perform two main steps to make this theory useful: first we need to introduce and study harmonic corrections of the solutions (\emph{i.e.}  new variables that differ up to harmonic functions to the original ones) and then we need to link the properties of these new variables to the old ones in a very particular manner to exploit the information carried out by (\ref{Formula_MesureMHD}), this second step will be achieved with a detailed study of the energy inequality satisfied by the Leray-type solutions of (\ref{EquationMHDoriginal}). Let us also stress here that, among other deep results -that will be highlighted in due time- it is the common framework of parabolic Morrey spaces developed in \cite{ChCHJ1} and \cite{ChCHJ2} that allows us to connect these two regularity theories.\\

Our main result which introduces a new notion of solutions for the MHD equations and generalizes the partial regularity theory is the following:
\begin{Theoreme}[Regularity for Dissipative solutions]\label{Theorem_main_original}
Let $\Omega$ be a bounded domain of the form given in (\ref{DefConjuntoOmega}) and $(\vU, P, \vB)$ be a weak solution on $\Omega$ of the MHD equations \eqref{EquationMHDoriginal}. Assume that
\begin{itemize}
\item[1)] we have that $(\vU, \vB, P, \vF, \vG)$ satisfies the conditions:
\begin{equation}\label{Hipotheses_Func}
\vU, \vB \in L^{\infty}_t L^{2}_x \cap L^{2}_t \dot{H}^{1}_x(\Omega), \quad \vF, \vG \in L^{2}_{t} H^1_x(\Omega), \quad P \in \mathcal{D}'(\Omega);
\end{equation}
\item[2)] the solution $(\vU, P, \vB)$ is dissipative, \emph{i.e.}, the quantity
\begin{eqnarray}
M&=&-\partial_t(|\vU|^2 + |\vB|^2 )+ \Delta (|\vU|^2 + |\vB|^2 )  - 2 (|\vn \otimes \vU|^2 + |\vn \otimes \vB|^2 - 2 \langle div (P \vU) \rangle \label{Formula_MesureDissipativeMHD}\\
& &- div \left( (|\vU|^2 + |\vB|^2) \vU \right) + 2 div ((\vU \cdot \vB) \vB ) + 2 (\vF \cdot \vU + \vG \cdot\vB),\notag
\end{eqnarray}
is well-defined as a distribution and is a locally finite non-negative measure on $\Omega$;
\item[3)] there exists a positive constant $\epsilon^{*}$ such that for some $\left(t_{0}, x_{0}\right) \in \Omega$, we have
\begin{equation}\label{Hypo_PetitesseGradMHD}
\limsup _{r \rightarrow 0} \frac{1}{r} \iint_{]t_{0}-r^{2}, t_{0}+r^{2}[ \times B\left(x_{0}, r\right)}|\vn \otimes \vU|^{2} + |\vn \otimes \vB|^{2} dxds <\epsilon^{*},
\end{equation}
\end{itemize}
then $(\vU, \vB)$ is locally bounded: for every parabolic ball $Q$ which is compactly supported in a small neighborhood of the point $(t_0, x_0)$, we have $(\vU, \vB) \in L^\infty_tL^\infty_x (Q)$.  And thus, following Remark \ref{Rem_RegulariteSerrin} we obtain that $\vU, \vB\in L^\infty_t\dot H^2_x\cap L^2_t\dot H^3_x$ over a small neighborhood of the point $(t_0, x_0)$.
\end{Theoreme}
Some important remarks are in order here. We first note that the assumption on $\vU, \vB$ and $\vF, \vG$ stated in the first point of the theorem are classical and can be obtained for any weak Leray-type solution of (\ref{EquationMHDoriginal}), however we only assume here that $P$ is a general distribution and this fact implies that we must define in a very specific way the quantity $M$ given in (\ref{Formula_MesureDissipativeMHD}): indeed we need to define the expression $\langle div(P \vU)\rangle$ as a very particular limit to make it meaningful (see formula \eqref{def_limlim} below). This general point of view about the pressure constitutes one of the main novelties of this article since it introduces a new class of solutions for the MHD equations, called \emph{dissipative solutions}, and it allows us to generalize all the previous results related to the partial regularity theory of such equations.  Let us mention that these type of dissipative solutions were studied for the classical Navier-Stokes equations in \cite{CML1} and \cite{CML}. To end these preliminaries remarks, note also that the condition (\ref{Hypo_PetitesseGradMHD}) is rather classical in the setting of the partial regularity theory (see Hypothesis (\ref{HypothesePetitesseGrad}) in Theorem \ref{Teo_CKNMHD}). Moreover, in the framework of the previous theorem, if we denote by $\Sigma_0$ the set of points for which we have the following behavior
$$\limsup _{r \rightarrow 0} \frac{1}{r} \iint_{]t_{0}-r^{2}, t_{0}+r^{2}[ \times B\left(x_{0}, r\right)}|\vn \otimes \vU|^{2} + |\vn \otimes \vB|^{2}   dx ds\geq \varepsilon^*,$$
then it can be shown, by a standard Vitali covering lemma, that the parabolic Hausdorff measure of the set $\Sigma_0$ is null, which means that this set of ``irregular points'' is actually very small (see Section 13.10 of the book \cite{PGLR1}). \\

\textbf{Outline of the paper.} In Section \ref{Section_Preliminaries} we fix some notation and we explain the main strategy that will be used in order to prove Theorem \ref{Theorem_main_original}: we will follow some steps and each one of the remaining sections will be devoted to the proof of these steps. Some technical lemmas are stated in Appendices. 
\section{Notation and strategy of the proof}\label{Section_Preliminaries}
We start this section by recalling some notation and useful facts. First we recall the notion of parabolic  H\"older and Morrey spaces and for this we consider the homogeneous space $(\mathbb{R}\times \R, d, \mu)$ where $d$ is the parabolic quasi-distance given by 
\begin{equation*}
d\big((t,x), (s,y)\big)=|t-s|^{\frac{1}{2}}+|x-y|,
\end{equation*}
and where $\mu$ is the usual Lebesgue measure $d\mu=dtdx$. Note that the homogeneous dimension is now $Q=5$. More details about the homogeneous spaces can be found in the books \cite{Folland}, \cite{PGLR1}, \cite{Triebel}. 
Associated to this distance, we  define homogeneous (parabolic) H\"older spaces $\dot{\mathcal{C}}^\alpha(\mathbb{R}\times \R, \R)$ where $\alpha\in ]0,1[$ by the following condition:
\begin{equation}\label{Holderparabolic}
\|\vphi\|_{\dot{\mathcal{C}}^\alpha}=\underset{(t,x)\neq (s,y)}{\sup}\frac{|\vphi(t,x)-\vphi(s,y)|}{\left(|t-s|^{\frac{1}{2}}+|x-y|\right)^\alpha}<+\infty,
\end{equation}
and this formula studies H\"older regularity in both time and space variables. Now, for $1< p\leq q<+\infty$, parabolic Morrey spaces $\mathcal{M}_{t,x}^{p,q}$ are defined as the set of measurable functions $\vphi:\mathbb{R}\times\R\longrightarrow \R$ that belong to the space $(L^p_{t,x})_{loc}$ such that $\|\vphi\|_{M_{t,x}^{p,q}}<+\infty$ where
\begin{equation}\label{DefMorreyparabolico}
\|\vphi\|_{\mathcal{M}_{t,x}^{p,q}}=\underset{x_{0}\in \R, t_{0}\in \mathbb{R}, r>0}{\sup}\left(\frac{1}{r^{5(1-\frac{p}{q})}}\int_{|t-t_{0}|<r^{2}}\int_{B(x_{0},r)}|\vphi(t,x)|^{p}dxdt\right)^{\frac{1}{p}}.
\end{equation}
These spaces are generalization of usual Lebesgue spaces, note in particular that we have $\mathcal{M}_{t,x}^{p,p}=L_{t,x}^p$. 
We refer the readers to the book \cite{Triebel1} for a general theory concerning the Morrey spaces and Hölder continuity.\\

To continue, we will introduce a change of variables that will allow us to work with a more symmetric expression of the equations (\ref{EquationMHDoriginal}). Indeed, following Elsasser \cite{elsasser1950}, we define
\begin{equation}\label{Def_ChangeVariable}
\vu= \vU + \vB,\quad \vb= \vU-\vB, \quad \vf= \vF + \vG\quad \mbox{and} \quad\vg = \vF - \vG,
\end{equation}
and then the original system \eqref{EquationMHDoriginal} becomes  
\begin{equation}\label{EquationMHD}
\begin{cases}
\partial_{t}\vu=\Delta \vu -(\vb\cdot\vn)\vu-\vn P+\vf,\quad div (\vu) = div (\vf)=0,\\[3mm]
\partial_{t}\vb=\Delta \vb -(\vu\cdot\vn)\vb- \vn P +\vg,\quad div (\vb)=div (\vg)=0,\\[3mm]
\vu(0,x)=\vu_{0}(x), div(\vu_0)=0,\quad \vb(0,x)=\vb_{0}(x), div(\vb_0)=0,
\end{cases}
\end{equation}
where, since $div(\vu)=div(\vb)=0$, we have that $P$ satisfies the equation 
\begin{equation}\label{EquationPression}
\Delta P= - \sum^3_{i,j= 1} \partial_i \partial_j (u_i b_j),
\end{equation}
and from this equation, we remark that the pressure $P$ is only determined by the couple  $(\vu, \vb)$. \\

Now, let $\Omega$ be a bounded subset of the form given in (\ref{DefConjuntoOmega}), we say that the couple $(\vu, \vb)\in L^\infty_tL^2_x\cap L^2_t\dot{H}^1_x(\Omega)$ satisfies the MHD equations (\ref{EquationMHD}) in the weak sense if for all $\vphi, \vpphi \in \mathcal{D}(\Omega)$ such that $div(\vphi)=div(\vpphi)=0$, we have
$$
\begin{cases}
\langle\partial_{t}\vu-\Delta \vu +(\vb\cdot\vn)\vu-\vf|\vphi \rangle_{\mathcal{D}'\times \mathcal{D}}=0,\\[3mm]
\langle\partial_{t}\vb-\Delta \vb +(\vu\cdot\vn)\vb-\vg|\vpphi \rangle_{\mathcal{D}'\times \mathcal{D}}=0,
\end{cases}
$$
note that if $(\vu, \vb)$ are solutions of the previous system, then due to the expression (\ref{EquationPression}) there exists a pressure $P$ such that (\ref{EquationMHD}) is fulfilled in $\mathcal{D}'$. We will work from now on with the following set of hypotheses for the functions $\vu, \vb, \vf,\vg$ and $P$:
\begin{equation}\label{Hypotheses_Travail}
\begin{split}
&\vu, \vb \in L^{\infty}_t L^{2}_x \cap L^{2}_t \dot{H}^{1}_x(\Omega), \quad \vf, \vg \in L^{2}_{t} H^1_x(\Omega), \quad P \in \mathcal{D}'(\Omega),\qquad \\[4mm]
&\limsup _{r \rightarrow 0} \frac{1}{r} \iint_{]t_{0}-r^{2}, t_{0}+r^{2}[ \times B\left(x_{0}, r\right)}|\vn \otimes \vu|^{2} + |\vn \otimes \vb|^{2} dxds <\epsilon^{*},\qquad (t_0, x_0)\in \Omega,
\end{split}
\end{equation}
which can easily be deduced from (\ref{Def_ChangeVariable}), (\ref{Hipotheses_Func}) and (\ref{Hypo_PetitesseGradMHD}) and where the point $(t_0, x_0)\in \Omega$ will be fixed from now on. In addition, for the Elsasser form MHD equations \eqref{EquationMHD}, the corresponding dissipative distribution turns out to be  
\begin{eqnarray}
\lambda&=&-\partial_t(|\vu|^2 + |\vb|^2 )+ \Delta (|\vu|^2 + |\vb|^2 )  - 2 (|\vn \otimes \vu|^2 + |\vn \otimes \vb|^2 ) - \langle div \big(P (\vu + \vb)\big)\rangle \label{Formula_MesureDissipativeMHD1}  \\  
& & - div ( |\vu|^2 \vb + |\vb|^2 \vu ) + 2 (\vf \cdot \vu + \vg \cdot\vb) \notag,
\end{eqnarray}
let us remark that the quantity above includes less nonlinear terms than \eqref{Formula_MesureDissipativeMHD} in Theorem \ref{Theorem_main_original} due to the symmetric property of the Elsasser formulation.\\

Once our framework is clear, we will explain now the strategy that will be implemented for proving Theorem \ref{Theorem_main_original}:

\begin{itemize}
\item{\bf Step 1} Using as a starting point the system (\ref{EquationMHD}), the variables $\vu, \vb, P$ and the initial data $\vf, \vg$ (every term under the hypotheses (\ref{Hypotheses_Travail}) above), we will derive from $\vu$ and $\vb$ two new variables, $\vv$ and $\vh$ respectively, that will be called \emph{harmonic corrections} of $\vu$ and $\vb$ as they differ (locally) from the original variables up to harmonic functions. Section \ref{Section_newvariable} will be devoted to a detailed study of these harmonic corrections since they are one key ingredient of our method: we will first show how the hypotheses on $\vu$ and $\vb$ are transmitted to the new variables $\vv$ and $\vh$ and then we will investigate the PDEs (which are of the same shape of the MHD equations) satisfied by these variables. We will see then that the global environment of the functions $\vv$ and $\vh$ is essentially better than the variables $\vu$ and $\vb$. 

\item{\bf Step 2} Next, we carry out a precise study of the energy inequalities satisfied by the variables $(\vu, P, \vb)$ and by the harmonic corrections $\vv$ and $\vh$. This step is crucial since it will allow us first to define in a very specific manner the distribution \textcolor{red}{(\ref{Formula_MesureDissipativeMHD1})} and then to transfer some information from the original system to the new one satisfied by $\vv$ and $\vh$. All these computations will be performed in Section \ref{Section_link}.
\item{\bf Step 3} The two previous steps allows us to obtain a better framework for the harmonic corrections $\vv$ and $\vh$. In this step we will see how to apply the usual partial regularity theory (as developed in  \cite{ChCHJ2}) in order to obtain an actual gain of regularity for these variables $\vv$ and $\vh$. Note that at this step, regularity in time and space  can still be obtained for the harmonic corrections $\vv$ and $\vh$.
\item{\bf Step 4} This fourth step is devoted to deduce from the previous step a gain of regularity in the original variables $\vu$ and $\vb$. 
\item{\bf Step 5} In this step, we will recover the regularity on the variables $\vU$ and $\vB$ and then Theorem \ref{Theorem_main_original} will be completely proven.  This will be done by using the local partial regularity obtained in \cite{ChCHJ1}. It is worth to remark here that by applying the local regularity Theorem \ref{Teo_SerrinMHD} we may loose some information in the time variable and we will only be able to obtain a gain regularity for the space variable.
\end{itemize}
The steps 3, 4 and 5 will be studied in Sections \ref{subsection_vh}, \ref{Sec_GainReguubb} and \ref{Sec_RegUB} respectively. 

\section{The harmonic corrections - Step 1}\label{Section_newvariable}
From the functions $\vu$ and $\vb$ we are going to derive two new variables $\vv$ and $\vh$ that locally (\emph{i.e.} over a small neighborhood of a point $(t_0,x_0)$ fixed in (\ref{Hypotheses_Travail}) above) differ from $\vu$ and $\vb$ only up to harmonic functions, but before we need to introduce some specific neighborhoods of $(t_0, x_0)$ that will fix our framework: for a radius $0<\rho$ small enough, let $Q_{\rho}$ be a bounded subset of $\Omega$ of the form
\begin{equation}\label{Def_ParaBallde}
Q_{\rho}:=]t_0-\rho^2,t_0+\rho^2[\times B(x_0,\rho).
\end{equation}
where $ B(x_0, \rho)$ is an open ball in $\R$ at center $x_0 \in \R$. When the context is clear, for usual (euclidean) balls, we will write $B_\rho$ instead of $B(x,\rho)$. We will also need a smaller subset $\Qro$ of the form (\ref{Def_ParaBallde}) above where $0< \rho_0 < \rho$. \\

Let us now construct a cut-off function $ \psi:\mathbb{R}\times \mathbb{R}^3\longrightarrow \mathbb{R}$ such that $\psi\in \mathcal{C}^{\infty}_0(\mathbb{R}\times \mathbb{R}^3, \mathbb{R})$ and 
\begin{equation}\label{cut_off}
 supp(\psi)\subset \Qr \quad \text{and} \quad \psi\equiv 1 \mbox{ on } \Qrone.
\end{equation} 
with $0<\rho_0 < \rho_1 < \rho$.  We can now define the harmonic corrections $\vv$ and $\vh$.
\begin{Definition} Assume that $\vu$ and $\vb$ satisfy the general hypotheses (\ref{Hypotheses_Travail}) and consider the localizing function $\psi$ given in (\ref{cut_off}). Then we define the variables $\vv$ and $\vh$ by:
\begin{equation}\label{new_variable}
\vv := - \frac{1}{\Delta} \vn \wedge (\psi \vn \wedge \vu ), \quad
\vh := - \frac{1}{\Delta} \vn \wedge (\psi \vn \wedge \vb ).
\end{equation}
\end{Definition}
Note that by the localizing properties of $\psi$ (in particular $\psi \equiv 1$ on $\Qro\subset Q_{\rho_1}$), by the vector identity $\vn \wedge ( \vn \wedge \vu ) = \vn (div \vu) - \Delta \vu$ and by the divergence free conditions for $\vu$ and $\vb$, we have the (local) identities
\begin{equation}\label{Identite_Harmoniques1}
\Delta \vv =  \Delta \vu, \quad \text{and}\quad 
\Delta \vh =  \Delta \vb, \quad \text{over} \quad \Qro,
\end{equation}
and these identities explain the denomination of \emph{harmonic corrections} given to the variables $\vv$ and $\vh$.\\

We state now some elementary facts on $\vv$ and  $\vh$ that can be easily deduced from the general hypotheses on $\vu$ and $\vb$.
\begin{Proposition}\label{Proposition_vh}
Assume that $\vu, \vb \in L^{\infty}_t L^{2}_x \cap L^{2}_t \dot{H}^{1}_x(\Omega)$, then the functions $\vv, \vh$ defined by the formula \eqref{new_variable} satisfy the following two facts: 
\begin{itemize}
\item[1)] the functions $\vv$ and $\vh$ are divergence free: $div(\vv) =0$ and $div(\vh) =0 $;
\item[2)]  we have $\vv, \vh \in L^{\infty}_t L^{2}_x \cap L^{2}_t H^{1}_x(\Qro)$.
\end{itemize}
\end{Proposition}
{\bf Proof.} For the first point we recall that the divergence of a curl is always null, then by definition of the functions $\vv$ and $\vh$ given in (\ref{new_variable}) the first point follows easily.\\

For the second point we will only show that $\vv \in L^{\infty}_t L^{2}_x (\Qro)$ and $\vv \in L^{2}_t H^{1}_x(\Qro)$ and we will omit the details for $\vh$ as the arguments follow the same lines. Thus, using the vector identity $ \psi \vn \wedge \vu=\vn \wedge (\psi \vu ) - \vn \psi \wedge \vu$ and (\ref{new_variable}) we can write 
\begin{equation}\label{Reecriture_vv1}
\vv =  - \frac{1}{\Delta} \vn \wedge \big(\vn \wedge (\psi \vu )\big) + \frac{1}{\Delta} \vn \wedge \big(\vn \psi \wedge \vu\big),
\end{equation}
and taking the $L^2_x$-norm of $\vv$ on the ball $\Bro$ we obtain that
\begin{eqnarray*}
\|\vv(t,\cdot)\|_{L^2_x (\Bro)} 
&\leq& \left\| \frac{1}{\Delta} \vn \wedge \big(\vn \wedge (\psi \vu )\big)(t,\cdot) \right\|_{L^2_x (\Bro)} + \left\| \frac{1}{\Delta} \vn \wedge \big(\vn \psi \wedge \vu\big)(t,\cdot)\right\|_{L^2_x (\Bro)}  \\
& \leq &\left\| \frac{1}{\Delta} \vn \wedge \big(\vn \wedge (\psi \vu )\big)(t,\cdot) \right\|_{L^2_x (\R)} + \left\| \frac{1}{\Delta} \vn \wedge \big(\vn \psi \wedge \vu\big)(t,\cdot)\right\|_{L^2_x (\R)} \\
& \leq &C\left(\|\psi \vu (t,\cdot)\|_{L^2_x (\R)} + \|\vn \psi \wedge \vu(t,\cdot) \|_{L^{\frac{6}{5}}_x (\R)}\right),
\end{eqnarray*}
where in the last estimate we used the Hardy–Littlewood–Sobolev inequality $\norm{(-\Delta)^{-\frac{1}{2}} f}_{L^2 (\R)} \leq C \norm{f}_{L^{\frac{6}{5}}(\R)}$. Now by the support and regularity properties of $\psi$ (see (\ref{cut_off})) and by Hölder's inequality we  obtain 
$$\|\vv(t,\cdot)\|_{L^2_x (\Bro)} 
\leq C \left(\|\psi (t,\cdot)\|_{L^\infty_x (\Br)} + \|\vn \psi(t,\cdot)\|_{L^{3}_x (\Br)} \right) \|\vu (t,\cdot)\|_{L^{2}_x (\Br)}.
$$
Then, taking the $L^\infty$ norm in the time variable, we get that
$$
\norm{\vv}_{L^\infty_tL^2_x (\Qro)} 
\leq C \left(\|\psi \|_{L^\infty_{t, x} (\Qr)} + \|\vn \psi\|_{L^\infty_t L^{3}_x (\Qr)} \right) \|\vu \|_{L^\infty_t L^{2}_x (\Omega)}\leq C_{\psi} \|\vu \|_{L^\infty_t L^{2}_x (\Omega)}<+\infty,
$$
and we can conclude that $\vv \in L^\infty_tL^2_x (\Qro)$.

Now, using again the identity (\ref{Reecriture_vv1}) for $\vv$ and taking the $H^1_x$-norm we have
\begin{eqnarray*}
\|\vv(t, \cdot)\|_{ H^1_x (\Bro)}&\leq &\left\| \frac{1}{\Delta} \vn \wedge \big(\vn \wedge (\psi \vu )\big) (t, \cdot)\right\|_{ H^1_x (\R)} + \left\| \frac{1}{\Delta} \vn \wedge \big(\vn \psi \wedge \vu\big)(t, \cdot)\right\|_{ H^1_x (\R)}\\
&\leq & C\|\psi \vu(t, \cdot)\|_{ H^1_x (\R)}+ \left\|  \frac{1}{\Delta} \vn \wedge \big(\vn \psi \wedge \vu\big)(t, \cdot)\right\|_{L^2_x (\R)} \\
& &+ \left\|\vn \otimes \left(\frac{1}{\Delta} \vn \wedge \big(\vn \psi \wedge \vu\big)\right)(t, \cdot)\right\|_{L^2_x (\R)} \\
&\leq & C\|(\psi \vu)(t, \cdot)\|_{ H^1_x (\R)}+ C\|  \vn \psi \wedge \vu(t, \cdot) \|_{L^{\frac{6}{5}}_x (\R)} + C\|\vn \psi \wedge \vu (t, \cdot)\|_{L^2_x (\R)},
\end{eqnarray*}
where we used in the last line the same Hardy–Littlewood–Sobolev inequality. Now by the support properties of the localizing function $\psi$ and by H\"older's inequality we obtain
\begin{eqnarray*}
\|\vv(t, \cdot)\|_{ H^1_x (\Bro)}&\leq &C\left(\|\psi(t,\cdot)\|_{L^\infty_x}+\|\vn \psi(t,\cdot)\|_{L^\infty_x}\right)\|\vu (t, \cdot)\|_{H^1_x (\Br)}\\
&+& C\left(\|\vn \psi(t, \cdot) \|_{L^3_x} + \|\vn \psi(t, \cdot)\|_{L^{\infty}_x} \right) \|\vu(t, \cdot) \|_{L^{2}_x (\Br)},
\end{eqnarray*}
thus, taking the $L^2_t$-norm in the time variable and observing that since $\vu \in L^\infty_t L^2_x\cap L^2_t \dot{H}^1_x (\Omega)$ one trivially gets $\vu \in L^2_t H^1_x (\Omega)$, we have
$$\|\vv\|_{L^2_t H^1_x (\Qro)} \leq C_{\psi} \|\vu \|_{L^2_t H^1_x (\Omega)}<+\infty,$$
and this concludes the proof of Proposition \ref{Proposition_vh}.\hfill$\blacksquare$\\

This first proposition is rather straightforward and we need to study in more detail the properties of the new variables $\vv$ and $\vh$ with respect to the original functions $\vu$ and $\vb$. To this end, we define the following quantities
\begin{equation}\label{Def_DifferenceHarmonique}
\vbe := \vu - \vv \qquad \mbox{and}\qquad \vga := \vb- \vh,
\end{equation}
we can see from the identity (\ref{Identite_Harmoniques1}) that $\vu$ is equal to $\vv$ on $\Qro$ up to the harmonic corrector $\vbe$ and thus we have $\Delta \vbe=0$ over $\Qro$ and we also have $\Delta \vga=0$ over $\Qro$. 
\begin{Remarque}\label{remark_bega}The restriction on the smaller set $\Qro$ is not essential as long as we are only interested in proving the previous Proposition \ref{Proposition_vh} and following the same ideas of Proposition \ref{Proposition_vh} we can show that $\vv, \vh \in L^{\infty}_t L^{2}_x \cap L^{2}_t H^{1}_x(\Qr)$, from which we deduce the following estimates:
$$\|\vbe\|_{L^\infty_tL^2_x (\Qr)} \leq C_{\psi} \|\vu \|_{L^\infty_tL^2_x (\Omega)}, \quad 	\|\vbe\|_{L^2_t H^1_x (\Qr)} \leq C_{\rho, \psi} \|\vu \|_{L^2_t H^1_x (\Omega)}$$
and
$$\|\vga\|_{L^\infty_tL^2_x (\Qr)} \leq C_{\psi} \|\vb \|_{L^\infty_tL^2_x (\Omega)}, \quad \|\vga\|_{L^2_t H^1_x (\Qr)} \leq C_{\rho, \psi} \|\vb \|_{L^2_t H^1_x (\Omega)}.	$$
Note that this restriction over the set $\Qro$ will be crucial in the following Proposition \ref{Lemma_correctors}.
\end{Remarque}
The next proposition states a stronger estimate for $\vbe$ and $\vga$.
\begin{Proposition}\label{Lemma_correctors}
Under the general hypotheses (\ref{Hypotheses_Travail}) over $\vu$ and $\vb$ and using the variables $\vv$ and $\vh$ given in (\ref{new_variable}), 
 then	the functions $\vbe $ and $\vga$ defined in (\ref{Def_DifferenceHarmonique}) satisfy $\vbe, \vga \in L^\infty_t Lip_x (\Qro)$.
\end{Proposition}
This is a very useful result as it gives information on the differences between the original variables $(\vu, \vb)$ and the harmonic corrections $(\vv, \vh)$.\\

\noindent{\bf Proof.} Using expression (\ref{Reecriture_vv1}) we can write 
$$\vbe =  \vu -\vv = \vu + \frac{1}{\Delta} \vn \wedge \left( \vn \wedge (\psi \vu ) - \vn \psi \wedge \vu \right),$$
since we have the identity $\vn \wedge \vn \wedge (\psi \vu)=\vn(\vu \cdot \vn \psi) - \Delta (\psi \vu)$, as $\vu$ is divergence free, we can write
$$\vbe = \vu+\frac{1}{\Delta} \left( \vn (\vu \cdot \vn \psi)\right)-\psi \vu - \frac{1}{\Delta} \vn \wedge (\vn \psi \wedge \vu ).$$
Now, since  $\psi \equiv 1$ on $\Qro$ we obtain over the set $\Qro$ the identity
\begin{equation}\label{betauv}
\vbe = \frac{1}{\Delta} \left( \vn (\vu \cdot \vn \psi)\right) - \frac{1}{\Delta} \vn \wedge (\vn \psi \wedge \vu).
\end{equation}
Thus, for all $(t,x)\in\Qro$, for any $i=1,2,3$ and denoting by $K$ the convolution kernel associated to the operator $\frac{1}{\Delta}$ (namely $\frac{1}{|x|}$), we have for each term of the formula (\ref{betauv}) above
\begin{equation*}
\left|\partial_{x_i}\left( \frac{1}{\Delta} \big( \vn (\vu \cdot \vn \psi)\big) \right)(t,x)\right|= \left|  \partial_{x_i} K * \big( \vn (\vu \cdot \vn \psi)\big) (t,x) \right| 
\leq \sum_{j, k=1}^3 \left|\partial_{x_i} \partial_{x_k}  K * \big(u_j \partial_j \psi \big) (t,x) \right|
\end{equation*}
and 
\begin{equation*}
\left| \partial_{x_i}\left( \frac{1}{\Delta} \vn \wedge (\vn \psi \wedge \vu)\right) (t, x)\right| 
= \left| \partial_{x_i}K * \left( \vn \wedge (\vn \psi \wedge \vu ) \right) (t, x)\right|
\leq \sum_{j,k,\ell=1}^{3} \left|\partial_{x_i} \partial_{x_k} K * ( \partial_j \psi u_\ell) (t, x)\right|.
\end{equation*}
Now we recall that by construction (see (\ref{cut_off})) we have $supp(\vn \psi) \subset B^c_{\rho_1}$ with $\rho_0 < \rho_1<\rho$, thus for $(t,x) \in \Qro$, we have for $i=1,2,3$, 
\begin{eqnarray}
|\partial_{x_i} \vbe(t,x)|
& \leq &\sum_{j, k=1}^3 \left|\partial_{x_i} \partial_{x_k}  K * \big(u_j \partial_j \psi \big) (t,x) \right| + \sum_{j, k, \ell=1}^{3} \left|\partial_{x_i} \partial_{x_k} K * ( \partial_j \psi u_\ell) (t, x)\right|\notag \\
&\leq &C\int_{\{|x-y| > \rho_1 - \rho_0 , y \in \Br\}} \frac{1}{|x-y|^3} |\vn \psi (t,y)| |\vu (t,y)| \,dy\notag \\
&\leq &C\norm{\vn \psi (t, \cdot)}_{L^2 (\Br)} \norm{\vu (t, \cdot)}_{L^2 (\Br)},\label{betau_Lip}
\end{eqnarray}
which implies that $\vbe \in L^\infty_t Lip_x (\Qro)$ by the properties of the localizing function $\psi$ and the fact that $\vu \in L^\infty_tL^2_x (\Omega)$. The estimate for $\vga \in L^\infty_t Lip_x (\Qro)$ can be deduced in a similar manner. This ends the proof of Proposition \ref{Lemma_correctors}.\hfill$\blacksquare$
\begin{Remarque}\label{remark_bega1}
Note that the conclusion of the Proposition \ref{Lemma_correctors} above can be obtained, \emph{mutatis mutandis}, in any proper subset of $\Qr$.
\end{Remarque}

Once we have obtained some preliminary information on the local behavior of the functions $\vv$ and $\vh$, we will study now the PDEs satisfied by them. In this sense we have the following proposition:
\begin{Proposition}\label{Proposition_newPressureForce}
Over the set $\Qro$ of the form (\ref{Def_ParaBallde}), the functions $\vv$ and $\vh$ defined by (\ref{new_variable}) satisfy the system:
\begin{equation}\label{EquationMHD_companion}
\begin{cases}
\partial_{t}\vv=\Delta \vv -(\vh\cdot\vn)\vv-\vn q+\vk,\\[3mm]
\partial_{t}\vh=\Delta \vh -(\vv\cdot\vn)\vh- \vn r +\vl.
\end{cases}
\end{equation}
where the functions $q, r, \vk$ and $\vl$ satisfy the following properties
\begin{enumerate}
\item[1)]  $q \in L^{\frac{3}{2}}_{t,x} (\Qro)$ and  $r \in L^{\frac{3}{2}}_{t,x}(\Qro)$,
\item[2)]  $div(\vk)=0$, $div(\vl)=0$ and $\vk -\vk_0 \in L^{2}_{t,x} (\Qro)$ and  $\vl - \vl_0 \in L^{2}_{t,x}(\Qro)$, where 
\begin{equation}\label{Def_k0_l0}
\vk_0=- \frac{1}{\Delta} \vn \wedge (\psi \vn \wedge \vf )\qquad \mbox{and}\qquad \vl_0:= - \frac{1}{\Delta} \vn \wedge (\psi \vn \wedge \vg ),
\end{equation} 
are the harmonic corrections associated to the external forces $\vf$ and $\vg$.
\end{enumerate}
\end{Proposition}
The system (\ref{EquationMHD_companion}) satisfied by the couple $(\vv, \vh)$ is very similar to the problem (\ref{EquationMHD}) satisfied by the couple $(\vu, \vb)$, however there are many deep differences between these two equations. Indeed:
\begin{itemize}
\item[$\bullet$] The main feature of the equations (\ref{EquationMHD_companion}) relies in the fact that the terms $\vn q$ and $\vn r$, which can be considered as ``modified pressures'', \emph{do not} depend on the initial pressure $P$. This fact is related to the definition of the new variables $\vv$ and $\vh$ given in (\ref{new_variable}): the first operation made over $\vu$ and $\vb$ is given the Curl operator $\vn \wedge$ that annihilates all the gradients. 
\item[$\bullet$] The ``modified pressures'' $q$ and $r$ in (\ref{EquationMHD_companion}) are obtained in a very particular way by using vector calculus identities and the main point here is that we can deduce some information on these objects (namely $q,r\in L^{\frac{3}{2}}_{t,x} $).
\item[$\bullet$] Although equations (\ref{EquationMHD_companion}) and (\ref{EquationMHD}) are very similar, the information available for them is completely different: in (\ref{EquationMHD_companion}) we do have some control over the ``modified pressures'' while in (\ref{EquationMHD}) the pressure is only a general distribution. 
\end{itemize}
{\bf Proof.} We begin by considering $\pat \vv$ over the set $\Qro$. Since $\pat \psi= 0$ over $\Qro$ and using the equation satisfied by $\vu$ (see (\ref{EquationMHD})) we have 
\begin{eqnarray}
\pat \vv&=&  \pat \left( - \frac{1}{\Delta} \vn \wedge (\psi \vn \wedge \vu )\right)= - \frac{1}{\Delta} \vn \wedge (\psi \;\pat(\vn   \wedge \vu ))\notag\\
&=&- \frac{1}{\Delta} \vn \wedge \left(\psi \left(\Delta (\vn \wedge \vu) -\vn \wedge ( (\vb\cdot\vn)\vu )+ \vn \wedge \vf \right) \right)\notag \\
& = &- \frac{1}{\Delta} \vn \wedge (\psi \Delta (\vn \wedge \vu))
+  \frac{1}{\Delta} \vn \wedge (\psi \vn \wedge ( (\vb\cdot\vn)\vu )  ) - \frac{1}{\Delta} \vn \wedge (\psi \vn \wedge \vf)\notag\\
&:= &I_1 + I_2 + \vk_0.\label{equationdet_timevv}
\end{eqnarray}
Doing the same procedure for $\pat \vh$, we get that 
\begin{eqnarray}
\pat \vh&=& - \frac{1}{\Delta} \vn \wedge (\psi \Delta (\vn \wedge \vb))
+  \frac{1}{\Delta} \vn \wedge (\psi \vn \wedge ( (\vu \cdot\vn)\vb )  ) - \frac{1}{\Delta} \vn \wedge (\psi \vn \wedge \vg)\notag\\
&:= &I_3 + I_4 + \vl_0. \label{equationdet_timevh}
\end{eqnarray}
Note that $div(\vk_0)=div(\vl_0)=0$ since the divergence of a curl is always null.\\

The equations (\ref{equationdet_timevv}) and (\ref{equationdet_timevh}) are still far away from the wished result (\ref{EquationMHD_companion}). In the following lines we will display some vectorial identities in order to force the terms $I_1, I_2$ and $I_3, I_4$ to be as in (\ref{EquationMHD_companion}). This process will of course bring up additional terms and some of them will be classified as \emph{pressures} an other as \emph{external forces}. Since the terms $I_1, I_2$ and $I_3,I_4$ are symmetric, we will only detail the computations for $I_1, I_2$.\\
\begin{itemize}
\item For the first term $I_1=- \frac{1}{\Delta} \vn \wedge (\psi \Delta (\vn \wedge \vu))$ given in (\ref{equationdet_timevv}), we use the classical identity 
$$\psi \Delta (\vn \wedge \vu)
= \Delta(\psi  (\vn \wedge \vu))+(\Delta \psi) (\vn \wedge \vu) -2\displaystyle{\sum_{i=1}^{3}\partial_{i}\left((\partial_{i}\psi) (\vn \wedge \vu)\right)},$$
and we obviously get that
\begin{eqnarray*}
-\frac{1}{\Delta} \vn \wedge (\psi \Delta (\vn \wedge \vu) )
&= &-\frac{1}{\Delta} \vn \wedge \left(\Delta(\psi \vn \wedge \vu)+(\Delta \psi)(\vn \wedge \vu) -2\displaystyle{\sum_{i=1}^{3}\partial_{i}\left((\partial_{i}\psi)(\vn \wedge \vu)\right)}\right)\\
&=& \Delta \left(-\frac{1}{\Delta} \vn \wedge (\psi \vn \wedge \vu)\right)-\frac{1}{\Delta} \vn \wedge \left((\Delta \psi)(\vn \wedge \vu) \right)\\
& &+2\frac{1}{\Delta} \vn \wedge \left(\sum_{i=1}^{3}\partial_{i}\left((\partial_{i}\psi)(\vn \wedge \vu)\right)\right)
\end{eqnarray*}
Recalling that $\vv=-\frac{1}{\Delta} \vn \wedge (\psi \vn \wedge \vu)$ (see expression (\ref{new_variable})) we can write
$$I_1=\Delta \vv-\frac{1}{\Delta} \vn \wedge \left((\Delta \psi)(\vn \wedge \vu) \right)+2\frac{1}{\Delta} \vn \wedge \left(\sum_{i=1}^{3}\partial_{i}\left((\partial_{i}\psi)(\vn \wedge \vu)\right)\right).$$
Now, applying the vector calculus identity 
\begin{equation}\label{identiteVectorielle}
\phi \vn \wedge \vec{Y} = \vn \wedge (\phi \vec{Y}) - (\vn \phi) \wedge \vec{Y},
\end{equation}
to the last two terms of the right-hand side of the expression of $I_1$ given above, we get that  
\begin{eqnarray*}
I_1&= &\Delta \vv-\frac{1}{\Delta} \vn \wedge\left(\vn \wedge((\Delta \psi) \vu)-\vn(\Delta \psi) \wedge \vu\right)\\
& &+2 \frac{1}{\Delta} \vn \wedge\left(\sum_{i=1}^{3} \partial_{i}\left(\vn \wedge\left(\left(\partial_{i} \psi\right) \vu\right) - (\vn(\partial_{i} \psi)) \wedge \vu \right)\right),
\end{eqnarray*}
which can be rewritten as
\begin{equation}\label{first_term_I1}
I_1 =\Delta \vv + \vk_1,
\end{equation}
with 
\begin{eqnarray}
k_1&=&-\frac{1}{\Delta} \vn \wedge\left(\vn \wedge((\Delta \psi) \vu)-\vn(\Delta \psi) \wedge \vu\right)\notag\\
& &+2 \frac{1}{\Delta} \vn \wedge\left(\sum_{i=1}^{3} \partial_{i}\left(\vn \wedge\left(\left(\partial_{i} \psi\right) \vu\right) - (\vn(\partial_{i} \psi)) \wedge \vu \right)\right).\label{expression_first_force}
\end{eqnarray}
Remark in particular that we have $div(\vk_1)=0$ as the divergence of a curl is always null.\\

\item We now treat the second term $I_2= \frac{1}{\Delta} \vn \wedge (\psi \vn \wedge ( (\vb\cdot\vn)\vu ))$ of (\ref{equationdet_timevv}).  Using again the vectorial identity (\ref{identiteVectorielle}), we can write
\begin{eqnarray*}
I_2& =& \frac{1}{\Delta} \vn \wedge \left(\vn \wedge (\psi ( (\vb\cdot\vn)\vu ))  - (\vn \psi ) \wedge ( (\vb\cdot\vn)\vu ) \right)\\
& = &\frac{1}{\Delta} \vn \wedge \vn \wedge (\psi ( (\vb\cdot\vn)\vu ))  
-  \frac{1}{\Delta} \vn \wedge \left((\vn \psi ) \wedge ( (\vb\cdot\vn)\vu ) \right)\\
&= &\frac{1}{\Delta} \left( \vn div (\psi ( (\vb\cdot\vn)\vu )) \right) 
- \psi  (\vb\cdot\vn)\vu -  \frac{1}{\Delta} \vn \wedge \left((\vn \psi ) \wedge ( (\vb\cdot\vn)\vu ) \right)
\end{eqnarray*}
where in the last line we used the identity $\vn \wedge \vn \wedge \vY = \vn div(\vY) - \Delta \vY$.  Let us now define $q_1$ and $\vk_2$ by the expressions
\begin{equation}\label{expression_k2}
q_1 = - \frac{1}{\Delta}  div (\psi ( (\vb\cdot\vn)\vu )) \qquad \mbox{and}\qquad
\vk_2 = -  \frac{1}{\Delta} \vn \wedge \left((\vn \psi ) \wedge ( (\vb\cdot\vn)\vu ) \right), 
\end{equation}
we can thus write:
\begin{equation}\label{second_term_I2}
I_2 = -\vn q_1 - \psi  (\vb\cdot\vn)\vu  + \vk_2,
\end{equation}
note here that we clearly have $div(\vk_2)= 0$ since again the divergence of a curl is always null.\\ 

We need now to treat the term $\psi (\vb\cdot\vn)\vu$ which is present in (\ref{second_term_I2}). Recalling the definitions of the functions $\vbe =\vu - \vv$ and $\vga= \vb- \vh$ given in (\ref{Def_DifferenceHarmonique}) and by the fact that $\psi \equiv 1$ on the set $\Qro$, then we have 
$$\psi \vu = \vv + \vbe, \quad \psi \vb = \vh + \vga \quad \text{on} \quad \Qro,$$
which implies that, on $\Qro$, we have the identities 
\begin{eqnarray*}
\psi (\vb\cdot\vn)\vu= ((\psi \vb ) \cdot\vn )( \psi \vu )
=  [(\vh + \vga) \cdot\vn](\vv + \vbe)
= (\vh  \cdot\vn ) \vv + \vA,
\end{eqnarray*}
where $ \vA := (\vh  \cdot\vn)\vbe +  (\vga \cdot\vn)\vv + 
(\vga \cdot\vn)\vbe $. Observe that on $\Qro$, we have $ \vA = \psi \vA $ and we can rewrite the quantity $\psi \vA$ in the following manner 
$$\psi \vA := \vn q_2 - \vk_3,$$
with 
$$q_2 =  \frac{1}{\Delta} div (\psi \vA) \quad \text{and}\quad \vk_3 =  -\psi \vA + \vn \frac{1}{\Delta} div (\psi \vA),$$
we then have 
\begin{equation}\label{Decomposition1}
\psi (\vb\cdot\vn)\vu = (\vh  \cdot\vn ) \vv + \vn q_2 - \vk_3.
\end{equation}
Notice also that $div(\vk_3)=0$, indeed 
$$div(\vk_3)=-div(\psi \vA) + div(\vn \frac{1}{\Delta} div (\psi \vA)))=-div(\psi \vA)+div(\psi \vA)=0.$$
Thus, substituting the expression (\ref{Decomposition1}) above into \eqref{second_term_I2} we obtain the following formula for the term $I_2$:
\begin{equation}\label{second_term_I2_2}
I_2 = -\vn q_1 - (\vh  \cdot\vn ) \vv - \vn q_2 + \vk_3 + \vk_2.\\[5mm]
\end{equation}
\end{itemize}
Now using the expression \eqref{first_term_I1} for $I_1$ and \eqref{second_term_I2_2} for $I_2$ and coming back to \eqref{equationdet_timevv} we finally obtain
\begin{equation*} 
\pat \vv
= \Delta \vv  - (\vh  \cdot\vn ) \vv -\vn (q_1+ q_2) + \vk_0 + \vk_1 + \vk_2 + \vk_3 .
\end{equation*}
Performing the same computations for the right-hand side of \eqref{equationdet_timevh}, we get a very similar equation for $\vh$ : 
\begin{equation*}
\pat \vh
= \Delta \vh  - (\vv  \cdot\vn ) \vh -\vn (r_1+ r_2) + \vl_0 + \vl_1 + \vl_2 + \vl_3
\end{equation*}
where
$$r_1 = - \frac{1}{\Delta}  div \left(\psi [ (\vu\cdot\vn)\vb ]\right), \quad r_2 =  \frac{1}{\Delta} div [\psi \left((\vv  \cdot\vn)\vga + (\vbe \cdot\vn)\vh + 
(\vbe \cdot\vn)\vga \right) ],$$ 
and
$$\vl_{1}=-\frac{1}{\Delta} \vn \wedge[\vn \wedge(\Delta \psi \vb)-\vn(\Delta \psi) \wedge \vb]+2 \frac{1}{\Delta} \vn \wedge\left[\sum_{j=1}^{3} \partial_{j}\left(\vn \wedge((\partial_{j} \psi) \vb ) - (\vn(\partial_{j} \psi)) \wedge \vb \right)\right],$$
$$\vl_2 = - \frac{1}{\Delta} \vn \wedge \left((\vn \psi ) \wedge ( (\vu\cdot\vn)\vb ) \right), \quad \vl_3 = \vn r_2 - \psi \left((\vv \cdot\vn)\vga + 
(\vbe \cdot\vn)\vh + (\vbe \cdot\vn)\vga\right).\\[5mm]$$
Up to now, we have obtained the following equations for the new couple $(\vv, \vh)$ over the set $\Qro$: 
\begin{equation*}
\begin{cases}
\partial_{t}\vv=\Delta \vv -(\vh\cdot\vn)\vv-\vn q+\vk,\\[3mm]
\partial_{t}\vh=\Delta \vh -(\vv\cdot\vn)\vh- \vn r +\vl,
\end{cases}
\end{equation*}
with $q= q_1 + q_2$, $r= r_1 + r_2$, $\vk= \vk_0 + \vk_1 +  \vk_2 + \vk_3$ and $\vl=  \vl_0 + \vl_1 + \vl_2 + \vl_3$ and where $div(\vk)=div(\vl)=0$.\\[5mm]

Once we have obtained the previous system for $\vv$ and $\vh$, we need to show that  we have the two points of the Proposition \ref{Proposition_newPressureForce}, \emph{i.e.} that $q, r \in L^{\frac 32}_{t,x}(\Qro)$ and $\vk -\vk_0 , \vl - \vl_0 \in L^{2}_{t,x}(\Qro)$. Let us start with the first point of Proposition \ref{Proposition_newPressureForce}, but we will only detail the estimates for $q$ and $\vk$, as the corresponding estimates for $r$ and $\vl$ follow the same lines. 

\begin{itemize}
\item[1)] Let us prove here that $q \in L^{\frac{3}{2}}_{t,x} (\Qro)$. Recall that $q=q_1+q_2$ with
$$q_1 =  - \frac{1}{\Delta}  div (\psi ( (\vb\cdot\vn)\vu )) \quad \text{and} \quad  q_2 = \frac{1}{\Delta} div [\psi \left((\vh  \cdot\vn)\vbe + 
(\vga \cdot\vn)\vv + (\vga \cdot\vn)\vbe\right)],$$
and we will study these two terms separately.\\

We start taking the $L^{\frac{3}{2}}_{x} (\Bro)$-norm of $q_1$, and we apply Hölder's inequality, the Hardy-Littlewood-Sobolev inequality $\norm{(-\Delta)^{-\frac{1}{2}} f}_{L^{\frac{9}{5}} (\R)} \leq C \norm{f}_{L^{\frac{9}{8}}(\R)}$ and we use the support properties of the function $\psi$ (recall (\ref{cut_off})) to obtain: 
\begin{eqnarray}
\norm{q_{1} (t, \cdot)}_{L^{\frac{3}{2}}_x (\Bro)}&\leq &C_{\rho_0} \|-\frac{1}{\Delta}  div (\psi ( (\vb\cdot\vn)\vu )) (t, \cdot) \|_{L^{\frac{9}{5}}_x (\Bro)} \leq C_{\rho_0} \|-\frac{1}{\Delta}  div (\psi ( (\vb\cdot\vn)\vu )) (t, \cdot) \|_{L^{\frac{9}{5}}_x (\mathbb{R}^3)}\notag \\
&\leq &C_{\rho_0} \|\psi ( (\vb\cdot\vn)\vu ) (t, \cdot) \|_{L^{\frac{9}{8}}_x (\mathbb{R}^3)} =C_{\rho_0} \|\psi ( (\vb\cdot\vn)\vu ) (t, \cdot) \|_{L^{\frac{9}{8}}_x (\Br)}\notag \\
&\leq &C_{\rho_0}\|\psi(t, \cdot)\|_{L^{\infty}_x (\Br)} \|\vb (t, \cdot) \|_{L^{\frac{18}{7}}_x (\Br)} \|\vu (t, \cdot) \|_{H^{1}_x (\Br)}.\label{estimate_q1}
\end{eqnarray}
By interpolation we have the estimate $\|\vb (t, \cdot)\|_{L^{\frac{18}{7}}_x (\Br)} \leq \|\vb (t, \cdot) \|^{\frac{2}{3}}_{L^{2}_x (\Br)}  \|\vb (t, \cdot) \|^{\frac{1}{3}}_{L^{6}_x (\Br)} $ and using the classical Sobolev embedding $H^1_x \hookrightarrow L^6_x$, we thus have 
\begin{equation*}
\begin{aligned}
\norm{q_{1} (t, \cdot)}_{L^{\frac{3}{2}}_x (\Bro)}
\leq C_{\rho_0} \|\psi(t, \cdot)\|_{L^{\infty}_x (\Br)}\|\vb (t, \cdot) \|^{\frac{2}{3}}_{L^{2}_x (\Br)}  \|\vb (t, \cdot) \|^{\frac{1}{3}}_{H^{1}_x (\Br)}  \|\vu (t, \cdot) \|_{H^{1}_x (\Br)}.
\end{aligned}
\end{equation*}
Now, taking the $L^{\frac{3}{2}}$-norm in the time variable and using Hölder's inequality, we finally obtain 
\begin{equation}\label{estimate_q11}
\norm{q_{1}}_{L^{\frac{3}{2}}_{t,x} (\Qro)}\leq C_{\rho_0, \psi} \|\vb \|^{\frac{2}{3}}_{L^\infty_t L^{2}_x (\Qr)}  \|\vb \|^{\frac{1}{3}}_{ L^2_t H^{1}_x (\Qr)} 
\|\vu \|_{L^2_t H^{1}_x (\Qr)}<+\infty.
\end{equation}

We study now the term $q_2$. We first take the $L^{\frac{3}{2}}_{x} (\Bro)$-norm of $q_2$ and following the same ideas displayed in the estimate \eqref{estimate_q1}, we obtain
\begin{eqnarray}
\norm{q_{2} (t, \cdot)}_{L^{\frac{3}{2}}_x (\Bro)}
&\leq& C_{\rho_0} \left\|\frac{1}{\Delta} div [\psi \left((\vh  \cdot\vn)\vbe + 
(\vga \cdot\vn)\vv +(\vga \cdot\vn)\vbe\right)] (t, \cdot) \right\|_{L^{\frac{9}{5}}_x (\Bro)} \notag\\
& \leq &C_{\rho_0} \left\|\psi \left((\vh  \cdot\vn)\vbe + (\vga \cdot\vn)\vv + (\vga \cdot\vn)\vbe\right) (t, \cdot) \right\|_{L^{\frac{9}{8}}_x (\Br)} \notag\\
&\leq &C_{\rho_0}\|\psi(t, \cdot)\|_{L^\infty_x(\Br)}\left(\|\vh (t, \cdot)\|_{L^{\frac{18}{7}}_x (\Br)} \|\vbe (t, \cdot) \|_{H^{1}_x (\Br)}\right.\notag\\
& &\qquad\left. +\|\vga (t, \cdot) \|_{L^{\frac{18}{7}}_x (\Br)} \big(\|\vv (t, \cdot) \|_{H^{1}_x (\Br)} +  \|\vbe (t, \cdot)\|_{H^{1}_x (\Br)} \big)\right).\label{estimate_q2}
\end{eqnarray}
Note now that by interpolation and Sobolev embedding, we have $\|\vh \|_{L^{\frac{18}{7}}_x (\Br)} 
\leq  \|\vh\|^{\frac{2}{3}}_{L^{2}_x (\Br)}  \|\vh \|^{\frac{1}{3}}_{H^{1}_x (\Br)} $ and using Remark \ref{remark_bega} we can write
$$\|\vga (t, \cdot) \|_{L^{\frac{18}{7}}_x (\Br)} \leq  \|\vga (t, \cdot)\|^{\frac{2}{3}}_{L^{2}_x (\Br)}  \|\vga (t, \cdot) \|^{\frac{1}{3}}_{H^{1}_x (\Br)} 
\leq \|\vb (t, \cdot) \|^{\frac{2}{3}}_{L^{2}_x (\Br)}  \|\vb (t, \cdot)\|^{\frac{1}{3}}_{H^{1}_x (\Br)}. $$
Since we also have $\norm{\vbe}_{H^1_x (\Br)} \leq C_{\rho, \psi} \|\vu \|_{H^1_x (\Br)}$, substituting all these three estimates into \eqref{estimate_q2}, one gets that
\begin{equation*}
\begin{aligned}
\norm{q_{2} (t, \cdot)}_{L^{\frac{3}{2}}_x (\Bro)}
\leq & C_{\rho_0, \psi} 
\|\vh (t, \cdot)\|^{\frac{2}{3}}_{L^{2}_x (\Br)}  \|\vh (t, \cdot) \|^{\frac{1}{3}}_{H^{1}_x (\Br)} \|\vu (t, \cdot) \|_{H^{1}_x (\Br)} \\
&+ C_{\rho_0, \psi} \left(\|\vb (t, \cdot) \|^{\frac{2}{3}}_{L^{2}_x (\Br)}  \|\vb (t, \cdot) \|^{\frac{1}{3}}_{H^{1}_x (\Br)} \left(\|\vv (t, \cdot) \|_{H^{1}_x (\Br)} + \|\vu (t, \cdot) \|_{H^{1}_x (\Br)}\right)\right).
\end{aligned}
\end{equation*}
Now, taking $L^{\frac{3}{2}}$-norm in the time variable, by H\"older's inequality, and recalling that by Remark \ref{remark_bega} we have $\vv, \vh \in L^{\infty}_t L^{2}_x \cap L^{2}_t H^{1}_x(\Qr)$, we can conclude that
\begin{eqnarray}
\norm{q_{2}}_{L^{\frac{3}{2}}_{t,x} (\Qro)}
&\leq & C_{\rho_0, \psi} \|\vh \|^{\frac{2}{3}}_{L^\infty_t L^{2}_x (\Qr)}  \|\vh \|^{\frac{1}{3}}_{ L^2_t H^{1}_x (\Qr)} \|\vu \|_{L^2_t H^{1}_x (\Qr)} \notag\\
& & + C_{\rho_0, \psi} \|\vb \|^{\frac{2}{3}}_{L^\infty_t L^{2}_x (\Qr)}  \|\vb \|^{\frac{1}{3}}_{ L^2_t H^{1}_x (\Qr)} \left(\|\vv \|_{L^2_t H^{1}_x (\Qr)} + \|\vu \|_{L^2_t H^{1}_x (\Qr)}\right)<+\infty. \label{estimate_q22}
\end{eqnarray} 
Since $q=q_1+q_2$, by the previous estimates \eqref{estimate_q11} and \eqref{estimate_q22}, we finally obtain that $q \in L^{\frac{3}{2}}_{t,x} (\Qro)$.\\

\item[2)] We study now the second point of Proposition \ref{Proposition_newPressureForce}. As mentionned above (precisely, the line below the system of $(\vv, \vh))$, we know that $div(\vk)=0$, so it only remains to prove that $\vk-\vk_0=\vk_1+\vk_2+\vk_3\in L^{2}_{t,x} (\Qro)$ and we will study the terms $\vk_1, \vk_2$ and $\vk_3$ separately.\\
\begin{itemize}
\item[$\bullet$]Let us now recall the expression of $\vk_1$ given in \eqref{expression_first_force}:
\begin{equation*}
\vk_{1}=-\frac{1}{\Delta} \vn \wedge[\vn \wedge(\Delta \psi \vu)-\vn(\Delta \psi) \wedge \vu]
+2 \frac{1}{\Delta} \vn \wedge\left[\sum_{i=1}^{3} \partial_{i}\left(\vn \wedge\left(\left(\partial_{i} \psi\right) \vu\right) - (\vn(\partial_{i} \psi)) \wedge \vu \right)\right].
\end{equation*}
Taking the $L^2_x$-norm of $\vk_1$ over $\Bro$, by the boundedness of Riesz transforms in Lebesgue spaces and by the Hardy-Littlewood-Sobolev inequalities, we obtain that 
\begin{eqnarray*}
\|\vk_{1} (t, \cdot)\|_{L^2_x (\Bro)}&\leq &  \left\|\frac{1}{\Delta} \vn \wedge\vn \wedge(\Delta \psi\,\vu)(t, \cdot)\right\|_{L^2_x (\mathbb{R}^3)}+\left\|\frac{1}{\Delta} \vn \wedge(\vn(\Delta \psi) \wedge \vu) (t, \cdot)\right\|_{L^2_x (\mathbb{R}^3)}\\
& &+2\sum_{i=1}^{3}\left\|  \frac{1}{\Delta} \vn \wedge \partial_{i}\left(\vn \wedge\left(\left(\partial_{i} \psi\right) \vu\right) - (\vn(\partial_{i} \psi)) \wedge \vu \right)(t, \cdot)\right\|_{L^2_x (\mathbb{R}^3)}\\
&\leq & C\left\|(\Delta \psi\, \vu)(t, \cdot)\right\|_{L^2_x (\mathbb{R}^3)}+C\left\|(\vn(\Delta \psi) \wedge \vu)(t, \cdot) \right\|_{L^{\frac65}_x (\mathbb{R}^3)}\\
& &+C\sum_{i=1}^{3}\left\|\left(\vn \wedge\left(\left(\partial_{i} \psi\right) \vu\right) - (\vn(\partial_{i} \psi)) \wedge \vu \right)(t, \cdot)\right\|_{L^2_x (\mathbb{R}^3)}.
\end{eqnarray*}
By the localization property of $\psi$ (see (\ref{cut_off})) and by H\"older's inequality we can write
\begin{eqnarray*}
\|\vk_{1} (t, \cdot)\|_{L^2_x (\Bro)}&\leq &C\|\Delta \psi(t, \cdot)\|_{L^\infty_x (\Br)}\left\|\vu(t, \cdot)\right\|_{L^2_x (\Br)}+C\|\vn(\Delta \psi)(t, \cdot)\|_{L^{3}_x (\Br)} \|\vu(t, \cdot) \|_{L^{2}_x (\Br)}\\
&&+C\|\Delta \psi(t, \cdot)\|_{L^\infty_x (\Br)}\|\vu(t, \cdot)\|_{L^2_x (\Br)}+C\|\vn \psi(t, \cdot)\|_{L^\infty_x (\Br)}\|\vu(t, \cdot)\|_{H^1_x (\Br)},
\end{eqnarray*}
thus, taking the $L^2$-norm in the time variable we obtain
\begin{equation}\label{estimate_k1}
\norm{\vk_{1}}_{L^2_{t, x} (\Qro)}\leq C_\psi \norm{\vu}_{L^2_t H^1_x (\Br)}<+\infty.
\end{equation}
\item[$\bullet$]For $\vk_2$, recalling its expression given in formula \eqref{expression_k2} and taking its $L^2_x$-norm on $\Bro$, we have by H\"older's inequality 
\begin{eqnarray*}
\norm{\vk_2 (t, \cdot)}_{L^2_x (\Bro)}&=& \left\|\frac{1}{\Delta} \vn \wedge \left((\vn \psi ) \wedge ( (\vb\cdot\vn)\vu ) \right) (t, \cdot) \right\|_{L^2_x (\Bro)} \\
&\leq &C_{\rho_0} \left \|\frac{1}{\Delta} \vn \wedge \left((\vn \psi ) \wedge ( (\vb\cdot\vn)\vu ) \right) (t, \cdot) \right\|_{L^\infty_x (\Bro)}
\end{eqnarray*}
We will use now some of the arguments displayed in the proof of Lemma \ref{Lemma_correctors}. Indeed, since $\vn \psi \subset B^c_{\rho_1}$ with $\rho_1 > \rho_0$, for $(t,x ) \in \Qro$ we have
\begin{equation*}
\begin{split}
\left| \frac{1}{\Delta} \vn \wedge \left((\vn \psi ) \wedge ( (\vb\cdot\vn)\vu ) \right)(t,x) \right|
&\leq  C\int_{\{|x-y| > \rho_1 - \rho_0 , y \in \Br\}} \frac{1}{|x-y|^2} |\vn \psi (t,y)| |(\vb \cdot \vn )\vu (t,y)| \,dy \\
&\leq C_{\rho_0, \rho_1} \norm{\vn \psi (t, \cdot)}_{L^\infty_x (\Br)} \norm{\vb (t, \cdot)}_{L^2_x (\Br)} \norm{\vu (t, \cdot)}_{H^1_x (\Br)}.
\end{split}
\end{equation*}
Now, taking the $L^2$-norm in the time variable of the previous inequality we obtain that
\begin{equation}\label{estimate_k2}
\norm{\vk_2}_{L^2_{t, x} (\Qro)} \leq C_{\rho_0, \rho_1, \psi}  \norm{\vb}_{L^\infty_t L^2_x (\Qr)} \norm{\vu}_{L^2_t H^1_x (\Qr)}<+\infty.
\end{equation}

\item[$\bullet$]Recall now the expression of $\vk_3$:
\begin{eqnarray}
\vk_3& = &- \psi \left((\vh  \cdot\vn)\vbe + (\vga \cdot\vn)\vv + (\vga \cdot\vn)\vbe\right) +\vn \left(\frac{1}{\Delta} div [\psi ((\vh  \cdot\vn)\vbe + (\vga \cdot\vn)\vv + (\vga \cdot\vn)\vbe ) ]\right) \notag \\
& = &- \psi \vA + \vn \frac{1}{\Delta} div (\psi \vA)\label{estimate_k31}
\end{eqnarray}
where we denoted $\vA : = (\vh  \cdot\vn)\vbe + (\vga \cdot\vn)\vv + (\vga \cdot\vn)\vbe$.
 
We treat the first term in the right-hand side of \eqref{estimate_k31}. Since $\psi \equiv 1$ on the set $\Qro$, we have the following estimates in the space variable:
\begin{eqnarray*}
\norm{\psi \vA(t, \cdot)}_{L^2_x (\Bro)}&=& \norm{\vA (t, \cdot)}_{L^2_x (\Bro)} =  \|\big((\vh  \cdot\vn)\vbe + (\vga \cdot\vn)\vv + 
(\vga \cdot\vn)\vbe \big)(t, \cdot)\|_{L^2_x (\Bro)}  \\
&\leq & \|\vh (t, \cdot)\|_{L^2_x (\Bro)} \| \vn \vbe (t, \cdot) \|_{L^\infty_x (\Bro)} +  \|\vga (t, \cdot)\|_{L^\infty_x (\Bro)} \| \vn \vv (t, \cdot)\|_{L^2_x (\Bro)} \\
& & + \|\vga (t, \cdot)\|_{L^2_x (\Bro)} \| \vn \vbe (t, \cdot) \|_{L^\infty_x (\Bro)}. 
\end{eqnarray*}
Then, taking the $L^2$-norm in the time variable and using Hölder's inequality we obtain 
\begin{eqnarray*}
\norm{\psi \vA}_{L^2_{t, x} (\Qro)} 
& \leq  &\|\vh \|_{L^2_{t,x} (\Qro)} \| \vbe \|_{L^\infty_t W^{1, \infty}_x (\Qro)} 
+  \|\vga \|_{L^\infty_{t,x} (\Qro)} \| \vv \|_{L^2_t H^1_x (\Qro)} \\
&& + \|\vga \|_{L^2_{t,x} (\Qro)} \|\vbe \|_{L^\infty_t W^{1, \infty}_x (\Qro)}.
\end{eqnarray*}
Now we recall that from Proposition \ref{Lemma_correctors} and Proposition \ref{Proposition_vh}, we have  $\vbe, \vga \in L^\infty_t Lip_x (\Qro)$ and $\vv, \vh \in L^{\infty}_t L^{2}_x \cap L^{2}_t H^{1}_x(\Qro)$ and these informations allow us to conclude that 
\begin{equation}\label{estimate_k32}
\norm{\psi \vA}_{L^2_{t, x} (\Qro)} < + \infty.
\end{equation}
It remains to study the second term in the right-hand side of \eqref{estimate_k31}.  For this, we introduce now a new cut-off function $\Phi \in \mathcal{D} (\mathbb{R}\times \R)$ such that 
$$supp(\Phi)\subset Q_{\rho_2} \quad \text{and} \quad \Phi\equiv 1 \mbox{ on } Q_{\rho_3},$$
where $0< \rho_0 < \rho_3 < \rho_2 < \rho_1 < \rho$ (recall the general definition of parabolic balls given in (\ref{Def_ParaBallde})). 
Since $supp(\psi)\subset \Qr$ and $\psi\equiv 1 \mbox{ on } \Qrone$ by (\ref{cut_off}), we see that $\psi \Phi \equiv 1$ on $Q_{\rho_3}$ and 
$$\ \text{supp} (\psi \Phi) \subset Q_{\rho_2}, \quad \text{and} \quad
\text{supp} (\psi (1-\Phi)) \subset Q_{\rho} \backslash Q_{\rho_3}.$$
With this cut-off function $\Phi$, we decompose the second term in the right-hand side of \eqref{estimate_k31} as 
\begin{equation}\label{estimate_k33}
\vn \frac{1}{\Delta} div (\psi \vA) =
\vn \frac{1}{\Delta} div (\psi \Phi \vA) + \vn \frac{1}{\Delta} div (\psi (1-\Phi) \vA) : = I_1 + I_2,
\end{equation}
and we treat the two terms separately. For $I_1$, taking its $L^2_{t,x}$-norm on $\Qro$ and applying the same strategy used for $\psi \vA$, we have 
\begin{eqnarray}
\left\|\vn \frac{1}{\Delta} div (\psi \Phi \vA)\right\|_{L^2_{t,x} (\Qro)} 
&\leq &C\norm{ \psi \Phi \vA }_{L^2_{t,x} (Q_{\rho_2})} 
= C\|\psi \Phi \big((\vh  \cdot\vn)\vbe +(\vga \cdot\vn)\vv + 
(\vga \cdot\vn)\vbe \big)\|_{L^2_{t,x} (Q_{\rho_2})} \notag \\
&\leq & \; C \left( \|\vh \|_{L^2_{t,x} (Q_{\rho_2})} \| \vbe \|_{L^\infty_t W^{1, \infty}_x (Q_{\rho_2})} + \|\vga \|_{L^\infty_{t,x} (Q_{\rho_2})} \| \vv \|_{L^2_t H^1_x (Q_{\rho_2})} 
\right.\notag\\
&& \left. + \|\vga \|_{L^2_{t,x} (Q_{\rho_2})} \|\vbe \|_{L^\infty_t W^{1, \infty}_x (Q_{\rho_2})}\right) < + \infty,\label{estimate_k34}
\end{eqnarray}
where we used, following Remark \ref{remark_bega1}, the fact that we have $\vbe, \vga \in L^\infty_t Lip_x (Q_{\rho_2})$.\\

The estimate of the term $I_2$ is more subtle and different arguments shall be used: indeed taking the $L^2$-norm in the space variable, we obtain 
\begin{equation*}
\left\|\vn \frac{1}{\Delta} div (\psi (1-\Phi) \vA) (t, \cdot) \right\|_{L^2_{x} (\Bro)} 
\leq C_{\rho_0} \norm{\vn \frac{1}{\Delta} div (\psi (1-\Phi) \vA) (t, \cdot) }_{L^\infty_{x} (\Bro)}.
\end{equation*}
Then, by the support property of $\psi (1- \Phi)$ for any point $(t,x) \in \Qro$ we have
$$ \left|\vn \frac{1}{\Delta} div (\psi (1-\Phi) \vA) (t,x)\right| \leq   \;  C\int_{\{|x-y| > \rho_3 - \rho_0 , y \in \Br\}} \frac{1}{|x-y|^3} |\psi (1-\Phi) (t,y)| |\vA(t,y)| \,dy,$$
thus, by a careful study of the support properties of each term we can deduce the following estimate
\begin{eqnarray*}
\left|\vn \frac{1}{\Delta} div (\psi (1-\Phi) \vA) (t,x)\right|&\leq &C_{\rho_0, \rho_3} \norm{\psi (1-\Phi) (t, \cdot)}_{L^\infty_x (\Br)}  \|\vA\|_{L^2_x (\Br)} \\
&\leq  & \; C_{\rho_0, \rho_3} \norm{\psi (1-\Phi) (t, \cdot)}_{L^\infty_x (\Br)} 
\left(\norm{\vh (t, \cdot)}_{L^2_x (\Br)} \norm{\vbe (t, \cdot)}_{H^1_x (\Br)} \right.\\
& &\left.
+ \norm{\vga (t, \cdot)}_{L^2_x (\Br)} \big( \norm{\vv (t, \cdot)}_{H^1_x (\Br)} +  \norm{\vbe (t, \cdot)}_{H^1_x (\Br)} \big) \right) 
\end{eqnarray*}
where we used Hölder's inequality for each term of $\vA= (\vh  \cdot\vn)\vbe + (\vga \cdot\vn)\vv + (\vga \cdot\vn)\vbe$. We can thus write
\begin{eqnarray*}
\left\|\vn \frac{1}{\Delta} div (\psi (1-\Phi) \vA) \right\|_{L^2_{t, x} (\Qro)}\notag 
&\leq  & \; C_{\rho_0, \rho_3, \psi, \Phi} \left(\norm{\vh}_{L^\infty_t L^2_x (\Qr)} \norm{\vbe }_{L^2_t H^1_x (\Qr)}\right. \notag\\
& + &\left.\norm{\vga }_{L^\infty_t L^2_x (\Qr)} \big( \norm{\vv }_{L^2_t H^1_x (\Qr)} +  \norm{\vbe}_{L^2_t H^1_x (\Qr)} \big) \right),
\end{eqnarray*}
since by Remark \ref{remark_bega} we have the inequalities $\norm{\vbe}_{L^2_t H^1_x (\Qr)} \leq C_{\rho, \psi} \|\vu \|_{L^2_t H^1_x (\Qr)}$ and 
$\norm{\vga}_{L^\infty_tL^2_x (\Qr)} \leq C_{\rho, \psi} \|\vb \|_{L^\infty_tL^2_x (\Qr)}$, then we obtain 
\begin{eqnarray}
\left\|\vn \frac{1}{\Delta} div (\psi (1-\Phi) \vA) \right\|_{L^2_{t, x} (\Qro)}&\leq & C_{\rho, \rho_0, \rho_3, \psi, \Phi} \left(\norm{\vh}_{L^\infty_t L^2_x (\Qr)} \norm{\vu }_{L^2_t H^1_x (\Qr)}\right. \notag\\
&+&\left. \norm{\vb }_{L^\infty_t L^2_x (\Qr)} \big( \norm{\vv }_{L^2_t H^1_x (\Qr)} +  \norm{\vu}_{L^2_t H^1_x (\Qr)} \big) \right)<+\infty.\label{estimate_k35}
\end{eqnarray}
With estimates \eqref{estimate_k34} and \eqref{estimate_k35} we obtain the term $\vn \frac{1}{\Delta} div (\psi \vA)$ given in (\ref{estimate_k33}) belongs to $L^{2}_{t,x} (\Qro)$ and with the estimate \eqref{estimate_k32} we obtain by (\ref{estimate_k31}) that $\vk_3\in L^{2}_{t,x} (\Qro)$.\\
\end{itemize}
\noindent With this previous information on $\vk_3$, and gathering together estimates \eqref{estimate_k1} and \eqref{estimate_k2}, we can finally conclude that 
$$\vk -\vk_0 \in L^{2}_{t,x} (\Qro).$$
\end{itemize}
This completes the proof of Proposition \ref{Proposition_newPressureForce}.
\hfill$\blacksquare$
\section{Link of original system and the companion equations - Step 2}\label{Section_link}
In this section, we aim to define the quantity \eqref{Formula_MesureDissipativeMHD1} by proving the following proposition. 
\begin{Proposition}\label{Prop_limlim}
Assume that $(\vu, P, \vb)$ is a weak solution on $Q_{\rho_0}$ of the MHD equations \eqref{EquationMHD} and that $(\vu, \vb, P, \vf, \vg)$ satisfies the conditions \eqref{Hypotheses_Travail}.\\
	
\noindent Let $\theta \in \mathcal{D} (\mathbb{R})$ and  $\varphi \in \mathcal{D} (\R)$ be two functions such that $\displaystyle{\int_{\mathbb{R}}}\theta(t) dt=1$, $\text{supp}  (\theta) \subset (-1, 1)$ and $\displaystyle{\int_{\R} \varphi(x) dx=1}$, $\text{supp}  (\varphi) \subset B(0, 1)$. We define the auxiliary function 
\begin{equation}\label{Mollifier}
\phi_{\alpha, \ep} = \theta_{\alpha} (t) \varphi_{\ep} (x),
\end{equation}
with the standard mollifiers $\theta_{\alpha}=  \frac{1}{\alpha} \theta(\frac{t}{\alpha})$, $\alpha>0$, and  $ \varphi_{\ep}=  \frac{1}{\ep^3} \varphi(\frac{x}{\ep})$, $\ep > 0$. Let $\widetilde{Q}_{\rho_0}$ be a subset of $Q_{\rho_0}$, then for $\alpha , \ep$ small enough, the distributions $\vu * \phi_{\alpha, \ep} $, $ \vb * \phi_{\alpha, \ep}$ and $P * \phi_{\alpha, \ep}$ are well-defined on $\widetilde{Q}_{\rho_0}$. Moreover, the limit
$$\lim_{\ep \to 0} \lim_{\alpha \to 0} div \left( P * \phi_{\alpha, \ep} (\vu * \phi_{\alpha, \ep} + \vb * \phi_{\alpha, \ep})  \right), $$
exists in $\mathcal{D}'(\widetilde{Q}_{\rho_0})$ and does not depend on the functions $\theta$ and $\varphi$.
\end{Proposition}
Assuming this proposition is true, we can introduce the  notation
\begin{eqnarray}\label{def_limlim}
\langle div \big(P (\vu + \vb)\big)\rangle := \lim_{\ep \to 0} \lim_{\alpha \to 0} div \left( P * \phi_{\alpha, \ep} (\vu * \phi_{\alpha, \ep} + \vb * \phi_{\alpha, \ep}) \right),
\end{eqnarray}
and thus, it is not hard to see that, within our framework, the quantity 
\begin{eqnarray*}
\lambda&=&-\partial_t(|\vu|^2 + |\vb|^2 )+ \Delta (|\vu|^2 + |\vb|^2 )  - 2 (|\vn \otimes \vu|^2 + |\vn \otimes \vb|^2 ) - div \left( |\vu|^2 \vb + |\vb|^2 \vu \right)\\  
& & - \langle div \big(P (\vu + \vb)\big)\rangle  + 2 (\vf \cdot \vu + \vg \cdot\vb),
\end{eqnarray*}
is a well-defined distribution which will lead us the following concept. 
\begin{Definition}[Dissipative solutions]\label{Def_dissi}
Let $(\vu, P, \vb)$ be a weak solution over some subset $\Omega$ of the form (\ref{DefConjuntoOmega}) of equations (\ref{EquationMHD}) that satisfy \eqref{Hypotheses_Travail}. We will say that $(\vu, P, \vb)$ is a \emph{dissipative} solution if the distribution $\lambda$ given above is a non-negative locally finite measure on $\Omega$.
\end{Definition}
We will see later on how to exploit this definition in order to prove the main results of this article. Let us prove now the proposition above.\\

\noindent{\bf Proof of Proposition \ref{Prop_limlim}.}
We begin by regularizing the MHD equations \eqref{EquationMHD} as follows: denoting 
$$\vu_{\alpha, \ep} = \vu * \phi_{\alpha, \ep},\; \vb_{\alpha, \ep} = \vb * \phi_{\alpha, \ep},\; P_{\alpha, \ep} = P * \phi_{\alpha, \ep},\; \vf_{\alpha, \ep} = \vf * \phi_{\alpha, \ep} \quad\mbox{and }\quad \vg_{\alpha, \ep} = \vg * \phi_{\alpha, \ep},$$
where $\phi_{\alpha, \ep}$ is the function given in (\ref{Mollifier}). Thus one has 
\begin{eqnarray*}
\partial_t(|\vu_{\alpha, \ep}|^2 + |\vb_{\alpha, \ep}|^2) 
&= &2 \vu_{\alpha, \ep} \cdot \pat \vu_{\alpha, \ep} + 2 \vb_{\alpha, \ep} \cdot \pat \vb_{\alpha, \ep}\\
&=&2 \vu_{\alpha, \ep} \cdot  \left(\Delta \vu -(\vb\cdot\vn)\vu-\vn P+\vf \right) * \phi_{\alpha, \ep} + 2 \vb_{\alpha, \ep} \cdot \left(\Delta \vb -(\vu\cdot\vn)\vb-\vn P+\vg \right)* \phi_{\alpha, \ep},
\end{eqnarray*}
which implies
\begin{equation}\label{regular_alep_ub}
\begin{aligned}
\partial_t(|\vu_{\alpha, \ep}|^2 + |\vb_{\alpha, \ep}|^2 )
&= \Delta (|\vu_{\alpha, \ep}|^2 + |\vb_{\alpha, \ep}|^2 )
- 2 (|\vn \otimes \vu_{\alpha, \ep}|^2 + |\vn \otimes \vb_{\alpha, \ep}|^2 )  - 2 div \left( P_{\alpha,\ep}  (\vu_{\alpha, \ep} + \vb_{\alpha, \ep}) \right)\\
& \hspace{-0.5cm} - 2 \vu_{\alpha, \ep} \cdot \left(  (\vb \cdot \vn  \vu ) * \varphi_{\alpha, \ep} \right) - 2 \vb_{\alpha, \ep} \cdot \left( (\vu \cdot \vn \vb ) * \varphi_{\alpha, \ep} \right) 
+ 2 (\vu_{\alpha, \ep} \cdot \vf_{\alpha, \ep} + \vb_{\alpha, \ep} \cdot \vg_{\alpha, \ep}),
\end{aligned}
\end{equation}
where we used divergence free conditions when dealing with the terms involving the pressure: indeed, we have $((\vn P) * \phi_{\alpha, \ep}) \cdot \vu_{\alpha, \ep}= div \left( (P * \phi_{\alpha, \ep}) \vu_{\alpha, \ep}\right) $ and $((\vn p) * \phi_{\alpha, \ep}) \cdot \vu_{\alpha, \ep}= div \left( (P * \phi_{\alpha, \ep}) \vb_{\alpha, \ep} \right) $.\\

By the same procedure used above in (\ref{regular_alep_ub}), we obtain the following regularized equation for the companion system \eqref{EquationMHD_companion}:
\begin{equation}\label{regular_alep_vh}
\begin{aligned}
\partial_t(|\vv_{\alpha, \ep}|^2 + |\vh_{\alpha, \ep}|^2 )
&= \Delta (|\vv_{\alpha, \ep}|^2 + |\vh_{\alpha, \ep}|^2 )
- 2 (|\vn \otimes \vv_{\alpha, \ep}|^2 + |\vn \otimes \vh_{\alpha, \ep}|^2 )  - 2 div \left( q_{\alpha,\ep}  \vv_{\alpha, \ep} + r_{\alpha,\ep}  \vh_{\alpha, \ep}\right)\\
& \hspace{-0.5cm} - 2 \vv_{\alpha, \ep} \cdot \left(  (\vh \cdot \vn  \vv ) * \varphi_{\alpha, \ep} \right) - 2 \vh_{\alpha, \ep} \cdot \left( (\vv \cdot \vn \vh ) * \varphi_{\alpha, \ep} \right) 
+ 2 (\vv_{\alpha, \ep} \cdot \vk_{\alpha, \ep} + \vh_{\alpha, \ep} \cdot \vl_{\alpha, \ep}).
\end{aligned}
\end{equation}
Now, we need some convergence lemmas to help us to pass to the limit $\alpha \to 0$. For the sake of simplicity, we use the notations
$$\vu_{\ep} = \vu * \varphi_{\ep},\quad \vb_{ \ep} = \vb * \varphi_{ \ep}, \quad
\vf_{ \ep} = \vf * \varphi_{ \ep},\quad \vg_{ \ep} = \vg * \varphi_{ \ep},$$
to indicate that only a mollification in the space variable is considered (see (\ref{Mollifier})).
\begin{Lemme}\label{lemme_limit_alpha}
We have the following strong convergences on the subset $Q_{\rho_0/4}$: 
\begin{itemize}
\item $\vu_{\alpha, \ep} \underset{\alpha \rightarrow 0}{\longrightarrow} \vu_{\ep}$ and $\vb_{\alpha, \ep} \underset{\alpha \rightarrow 0}{\longrightarrow} \vb_{\ep}$ in $L_{t,x}^{2}(Q_{\rho_0/4}) \bigcap L_{t}^{2} H_{x}^{1}(Q_{\rho_0/4})$,
\item $\vec{\nabla} \otimes \vu_{\alpha, \ep} \underset{\alpha \rightarrow 0}{\longrightarrow} \vn \otimes \vu_{\ep}$ and $\vn \otimes \vb_{\alpha, \ep} \underset{\alpha \rightarrow 0}{\longrightarrow} \vec{\nabla} \otimes \vb_{\ep}$ in $L_{t,x}^{2}(Q_{\rho_0/4})$,
\item $(\vu \otimes \vb) * \phi_{\alpha, \ep} \underset{\alpha \rightarrow 0}{\longrightarrow}(\vu \otimes \vb) * \varphi_{\ep}$ and $(\vb \otimes \vu) * \phi_{\alpha, \ep} \underset{\alpha \rightarrow 0}{\longrightarrow}(\vb \otimes \vu) * \varphi_{\ep}$ in $L_{t,x}^{2} (Q_{\rho_0/4})$,
\item  $\vf_{\alpha, \ep} \underset{\alpha \rightarrow 0}{\longrightarrow} \vf_{\ep}$ and $\vg_{\alpha, \ep} \underset{\alpha \rightarrow 0}{\longrightarrow} \vg_{\ep}$ in $L_{t,x}^{2} (Q_{\rho_0/4})$.
\end{itemize}
\end{Lemme}
We omit the proof here, since one can easily check that these convergences hold by the definition of functions $\phi_{\alpha, \ep}$ and $\varphi_\ep$ given in (\ref{Mollifier}) and by the properties of $\vu, \vb, \vf, \vg$.\\ 

Note that a completelty similar result can be obtained for the variables $\vv, \vh$ and the forces $\vk, \vl$, which allows us to pass to the limit $\alpha \to 0$ for the companion equations. 
\begin{Remarque}\label{RemarqueStrongCV}
Observe however, that by the first point of Proposition \ref{Proposition_newPressureForce}, we have an additional strong convergence for the terms $q$ and $r$, indeed:
\begin{equation}\label{conver_alpha}
\lim\limits_{\alpha \to 0} q * \phi_{\alpha, \ep} = q * \varphi_{\ep}, \quad 
\lim\limits_{\alpha \to 0} r * \phi_{\alpha, \ep} = r * \varphi_{\ep}, \quad \text{in} \quad L_{t,x}^{\frac{3}{2}} (Q_{\rho_0/4}).
\end{equation}
\end{Remarque}
Now, with the help of Lemma \ref{lemme_limit_alpha}, we are able to pass to the limit $\alpha \to 0$ for all the terms in equality \eqref{regular_alep_ub}, except for the term involving pressure. Namely,
\begin{equation}\label{regular_ep_ub}
\begin{aligned}
\partial_t(|\vu_{\ep}|^2 + |\vb_{\ep}|^2 )
&= \Delta (|\vu_{\ep}|^2 + |\vb_{\ep}|^2 )
- 2 (|\vn \otimes \vu_{\ep}|^2 + |\vn \otimes \vb_{\ep}|^2 ) 
- 2 \vu_{ \ep} \cdot \left(  (\vb \cdot \vn  \vu ) * \varphi_{\ep} \right)  \\  
& \hspace{0.5cm}
- 2 \vb_{\ep} \cdot \left( (\vu \cdot \vn \vb ) * \varphi_{\ep} \right) - 2 \lim_{\alpha \to 0} div \left( P_{\alpha,\ep}  (\vu_{\alpha, \ep} + \vb_{\alpha, \ep}) \right)  + 2 (\vu_{ \ep} \cdot \vf_\ep+ \vb_{\ep} \cdot \vg_\ep).
\end{aligned}
\end{equation}
We can also pass to the limit $\alpha \to 0$ for every term involving in \eqref{regular_alep_vh}, 
\begin{equation}\label{regular_ep_vh}
\begin{aligned}
\partial_t(|\vv_{\ep}|^2 + |\vh_{\ep}|^2 )
&= \Delta (|\vv_{\ep}|^2 + |\vh_{\ep}|^2 )
- 2 (|\vn \otimes \vv_{\ep}|^2 + |\vn \otimes \vh_{\ep}|^2 ) 
- 2 \vv_{ \ep} \cdot \left(  (\vh \cdot \vn  \vv ) * \varphi_{\ep} \right) \\  
& \hspace{0.5cm}
 - 2 \vh_{\ep} \cdot \left( (\vv \cdot \vn \vh ) * \varphi_{\ep} \right)  - 2div \left( (q * \varphi_{\ep}) \vv_{\ep} + (r * \varphi_{\ep}) \vh_{\ep} \right)  + 2 (\vv_{\ep} \cdot \vk_\ep+ \vh_{\ep} \cdot \vl_\ep)
\end{aligned}
\end{equation}
Notice that, since we just assume that $P \in \mathcal{D}'$, the limit term $\lim\limits_{\alpha \to 0} div \left( P_{\alpha,\ep}  (\vu_{\alpha, \ep} + \vb_{\alpha, \ep}) \right)$ involving the pressure remains in (\ref{regular_ep_ub}). Remark also that thanks to the additional strong convergence \eqref{conver_alpha}, there is no such term in \eqref{regular_ep_vh}.\\

We study now the limit when $\ep \to 0$ with the following lemma.
\begin{Lemme}\label{lemme_limit_ep}
We have the following strong convergence: 
\begin{itemize}
\item $\vu_{\ep} \underset{\ep \rightarrow 0}{\longrightarrow} \vu$ and $\vb_{\ep} \underset{\ep \rightarrow 0}{\longrightarrow} \vb$ in $L_{t,x}^{2} (Q_{\rho_0/4}) \bigcap L_{t}^{2} \dot H_{x}^{1}(Q_{\rho_0/4})$,
\item $\vec{\nabla} \otimes \vu_{\ep} \underset{\ep \rightarrow 0}{\longrightarrow} \vn \otimes \vu$ and $\vn \otimes \vb_{\ep} \underset{\ep \rightarrow 0}{\longrightarrow} \vec{\nabla} \otimes \vb$ in $L_{t,x}^{2} (Q_{\rho_0/4})$,
\item  $\vf_{\ep} \underset{\ep \rightarrow 0}{\longrightarrow} \vf$, $\vg_{\ep} \underset{\ep \rightarrow 0}{\longrightarrow} \vg$ in $L_{t,x}^{2} (Q_{\rho_0/4})$.
\end{itemize}
\end{Lemme}
The proof of this lemma is again straightforward and we omit the details.\\ 

As pointed out by the Remark \ref{RemarqueStrongCV} and by the first point of Proposition \ref{Proposition_newPressureForce}, for new pressure terms $q$ and $r$, we have the following strong convergence:  
 \begin{equation}\label{conver_ep}
 \lim\limits_{\ep \to 0} q * \varphi_{\ep} = q, \quad 
 \lim\limits_{\ep \to 0} r * \varphi_{\ep} = r, \quad \text{in} \quad L_{t,x}^{\frac{3}{2}} (Q_{\rho_0/4})
 \end{equation}
which allows us to pass to the limit for the term involving $q$ and $r$ in equality \eqref{regular_ep_vh}.\\ 

In order to deal with the general limit when $\ep,\alpha\to 0$, we introduce the following two distributions: 
\begin{equation}\label{Def_MuEspNuEsp}
\begin{split}
\mu_\ep &=  2 \vu_{ \ep} \cdot \left(  (\vb \cdot \vn  \vu ) * \varphi_{\ep} \right) + 2 \vb_{\ep} \cdot \left( (\vu \cdot \vn \vb ) * \varphi_{\ep} \right)  
 -  div (|\vu|^2 \vb + |\vb|^2 \vu)\\
&and\\
\eta_\ep &=  2 \vv_{ \ep} \cdot \left(  (\vh \cdot \vn  \vv ) * \varphi_{\ep} \right) + 2 \vh_{\ep} \cdot \left( (\vv \cdot \vn \vh ) * \varphi_{\ep} \right)  
- div (|\vv|^2 \vh + |\vh|^2 \vv).
\end{split}
\end{equation}
With these two quantities at hand, we pass to the limit $\ep \to 0$ for the equalities \eqref{regular_ep_ub} and \eqref{regular_ep_vh} to get 
\begin{equation}\label{regular_ub}
\begin{aligned}
\partial_t(|\vu|^2 + |\vb|^2 )
&= \Delta (|\vu|^2 + |\vb|^2 )
- 2 (|\vn \otimes \vu|^2 + |\vn \otimes \vb|^2 ) - div (|\vu|^2 \vb + |\vb|^2 \vu) \\ 
& \hspace{0.5cm} - 2 \lim_{\ep \to 0} \lim_{\alpha \to 0} div \left( P_{\alpha,\ep}  (\vu_{\alpha, \ep} + \vb_{\alpha, \ep}) \right)  - \lim_{\ep \to 0} \mu_\ep
 + 2 (\vu \cdot \vf+ \vb \cdot \vg),
\end{aligned}
\end{equation}
and
\begin{equation}\label{regular_vh}
\begin{aligned}
\partial_t(|\vv|^2 + |\vh|^2 )
&= \Delta (|\vv|^2 + |\vh|^2 )
- 2 (|\vn \otimes \vv|^2 + |\vn \otimes \vh|^2 ) - div (|\vv|^2 \vh + |\vh|^2 \vv) \\ 
& \hspace{0.5cm} -2  div \left( q \vv + r \vh \right) - \lim_{\ep \to 0} \eta_\ep
+ 2 (\vv \cdot \vk+ \vh \cdot \vl)  ,
\end{aligned}
\end{equation}

Recall that we aim to give a sense to the quantity 
$$
\lim_{\ep \to 0} \lim_{\alpha \to 0} div \left(  P_{\alpha,\ep}  (\vu_{\alpha, \ep} + \vb_{\alpha, \ep})  \right).
$$

In order to achieve this goal, we shall need the following proposition of $\mu_\ep, \eta_\ep$, which links essentially the original MHD system and the companion system. 
\begin{Proposition}\label{lemme_limit}
We have the following convergence in $\mathcal{D}' (Q_{\rho_0/4})$ for $\mu_\ep$ and $\eta_\ep$: 
\begin{equation}\label{limit_muetq}
\lim_{\ep \to 0} \mu_\ep - \eta_\ep = 0. 
\end{equation}
\end{Proposition}
Before giving the detailed proof of Proposition \ref{lemme_limit}, let us first assume this lemma is proven and let us continue the proof of Proposition \ref{Prop_limlim}.\\

Without loss of generality, we set $\widetilde{Q}_{\rho_0} := Q_{\rho_0/4}$, then $\vu * \phi_{\alpha, \ep} $, $ \vb * \phi_{\alpha, \ep}$ and $P * \phi_{\alpha, \ep}$ are well-defined on $\widetilde{Q}_{\rho_0}$. Indeed, if we choose $\ep < \frac{\rho_0}{2}$, then for any point $(t,x) \in Q_{\rho_0/4}$, we have $Q_{t,x,\ep} = ]t-\ep^2, t+\ep^2[ \times B(x, \ep)\subset Q_{\rho}$.\\

Next, we rewrite equalities \eqref{regular_ub} and \eqref{regular_vh} in the following manner:
\begin{eqnarray*}
\lim_{\ep \to 0} \mu_\ep +2 \lim_{\ep \to 0} \lim_{\alpha \to 0} div \left(  P_{\alpha,\ep}  (\vu_{\alpha, \ep} + \vb_{\alpha, \ep})\right)&=& -\partial_t(|\vu|^2 + |\vb|^2 ) + \Delta (|\vu|^2 + |\vb|^2 )\\
& &- 2 (|\vn \otimes \vu|^2 + |\vn \otimes \vb|^2 )- div (|\vu|^2 \vb + |\vb|^2 \vu)+ 2 (\vu \cdot \vf+ \vb \cdot \vg),  
\end{eqnarray*}
and 
\begin{eqnarray*}
- \lim_{\ep \to 0} \eta_\ep &=& \partial_t(|\vv|^2 + |\vh|^2 )- \Delta (|\vv|^2 + |\vh|^2 )+ 2 (|\vn \otimes \vv|^2 + |\vn \otimes \vh|^2 ) + div (|\vv|^2 \vh + |\vh|^2 \vv) \\
&&+2  div \left( q \vv + r \vh \right) -2 (\vv \cdot \vk+ \vh \cdot \vl).
\end{eqnarray*}
Thus summing these two terms we can write
\begin{eqnarray*}
\lim_{\ep \to 0} (\mu_\ep-\eta_\ep) +2 \lim_{\ep \to 0} \lim_{\alpha \to 0} div \left(  P_{\alpha,\ep}  (\vu_{\alpha, \ep} + \vb_{\alpha, \ep})\right)&=& -\partial_t(|\vu|^2 + |\vb|^2 ) + \Delta (|\vu|^2 + |\vb|^2 )\\
& &- 2 (|\vn \otimes \vu|^2 + |\vn \otimes \vb|^2 )- div (|\vu|^2 \vb + |\vb|^2 \vu)\\
&&+ 2 (\vu \cdot \vf+ \vb \cdot \vg)+\partial_t(|\vv|^2 + |\vh|^2 )- \Delta (|\vv|^2 + |\vh|^2 ) \\
&&+ 2 (|\vn \otimes \vv|^2 + |\vn \otimes \vh|^2 ) + div (|\vv|^2 \vh + |\vh|^2 \vv)\\
&&+2  div \left( q \vv + r \vh \right) -2 (\vv \cdot \vk+ \vh \cdot \vl). 
\end{eqnarray*}
Now we use the fact stated in formula \eqref{limit_muetq}, which is the conclusion of Proposition \ref{lemme_limit}, and we obtain the following identity
\begin{eqnarray*}
2 \lim_{\ep \to 0} \lim_{\alpha \to 0} div \left(  P_{\alpha,\ep}  (\vu_{\alpha, \ep} + \vb_{\alpha, \ep})\right)
&=& -\partial_t(|\vu|^2 + |\vb|^2 ) + \Delta (|\vu|^2 + |\vb|^2 )\\
& &- 2 (|\vn \otimes \vu|^2 + |\vn \otimes \vb|^2 )- div (|\vu|^2 \vb + |\vb|^2 \vu)\\
&&+ 2 (\vu \cdot \vf+ \vb \cdot \vg)+\partial_t(|\vv|^2 + |\vh|^2 )- \Delta (|\vv|^2 + |\vh|^2 ) \\
&&+ 2 (|\vn \otimes \vv|^2 + |\vn \otimes \vh|^2 ) + div (|\vv|^2 \vh + |\vh|^2 \vv)\\
&&+2  div \left( q \vv + r \vh \right) -2 (\vv \cdot \vk+ \vh \cdot \vl),\end{eqnarray*}
which shows the existence of the limit $ \lim\limits_{\ep \to 0} \lim\limits_{\alpha \to 0} div \left(  P_{\alpha,\ep}  (\vu_{\alpha, \ep} + \vb_{\alpha, \ep})\right)$ in $\mathcal{D}'(Q_{\rho_0/4})$ since all the terms present in the right-hand side of the previous formula are well defined and this ends the proof of Proposition \ref{Prop_limlim}.
\hfill$\blacksquare$\\

The rest of this section is devoted to the proof of Proposition \ref{lemme_limit} which mainly relies on an idea of Duchon and Robert in \cite{DuchonRobert} which studied some deep properties of energy dissipation for weak solutions in the Navier-Stokes equations. Let us remark that a similar result for the MHD equations has been given in \cite{GaoTan}, however,  for the sake of completeness, we state and prove the corresponding results in the Elsasser formulation of MHD system \eqref{EquationMHD}. \\

\noindent{\bf Proof of Proposition \ref{lemme_limit}.} In order to simply the notations, we introduce the following operator, which denotes the difference of the translation in space variable of a function $\vX : \mathbb{R} \times \R \longrightarrow \R$ and the function $\vX$ itself: 
$$
\delta_{y}[\vX ](t, x) : = \vX (t, x-y) - \vX (t,x).
$$
We will need the following general result:
\begin{Lemme}\label{prop_NRST}
Let $\vX, \vY, \vZ \in L^{3}_{t,x} (Q_{\rho_0})$ such that $div(\vX)=div(\vY)=div( \vZ)=0$.  Let $\theta \in \mathcal{D}(\R)$ be a smooth function on $\R$ such that $\displaystyle{\int_{\R}} \theta(x)dx =1$ and $supp(\theta) \subset B(0,1)$. We define $\theta_\ep (x)= \frac{1}{\ep^3} \theta (\frac{x}{\ep})$ with $\ep > 0$.
Then we define the following four quantities for all $(t,x) \in Q_{\rho_0/4}  \subset Q_{\rho_0}$:  
\begin{equation}\label{four_quantities}
\begin{aligned}
&N_{\ep} (\vX, \vY, \vZ) (t, x)=(\vY * \theta_{\ep}) \cdot\left([(\vX \cdot \vn ) \vZ] * \theta_{\ep}\right) + (\vZ * \theta_{\ep}) \cdot\left([(\vX \cdot \vn ) \vY] * \theta_{\ep}\right)-div((\vY \cdot \vZ) \vX)\\
&R_{\ep} (\vX, \vY, \vZ) (t, x)=\int_{\R}  \left(\vn  \theta_{\ep}(y) \cdot \delta_{y}[\vX ](t, x)  \right) \left(\vY(x-y)-\vZ(x)\right) \cdot \left(\vZ (x-y) - \vY (x)\right)d y\\
&S_{\ep} (\vX, \vY, \vZ) (t, x)= \int_{\R} \left( \vn \theta_{\ep}(y) \cdot \delta_{y}[\vX ](t, x) \right) \left(\delta_{y}[\vY ](t, x)\right) \cdot\left((\vZ * \theta_\ep)(x)-\vY(x)\right) d y\\
&T_{\ep} (\vX, \vY, \vZ)  (t, x)= \int_{\R} \left(\vn
\theta_{\ep}(y) \cdot \delta_{y}[\vX ](t, x)\right) \left(\delta_{y}[\vZ ](t, x) \right) \cdot\left((\vY * \theta_\ep)(x)-\vZ(x)\right) d y
\end{aligned}
\end{equation}
with $0<\ep< \frac{\rho_0}{2}$.
Then, we have the limit 
$$\lim_{\ep \to 0} (N_{\ep}+R_{\ep}-S_{\ep}-T_{\ep}) (\vX, \vY, \vZ)= 0.$$
\end{Lemme}
The proof of Lemma \ref{prop_NRST} is technical and can be found in the Appendix below. The previous lemma will not be enough for our purposes and the following result gives a more precise limit of the previous quantities $(S_{\ep} + T_{\ep} - R_{\ep}) (\vX, \vY, \vZ)$ under an additional condition on $\vX, \vY, \vZ$.
\begin{Lemme}\label{lemma_Tep}
Let $\vX, \vY, \vZ \in L_{t,x}^3 (\Qro)$ with $div(\vX)=div(\vY)=div( \vZ)=0$. Assume that at least one of them belongs to $L_t^\infty Lip_x (\Qro)$, then we have 
$$\lim_{\ep \to 0}(S_{\ep} + T_{\ep} - R_{\ep}) (\vX, \vY, \vZ) = 0 \quad \text{in} \quad L^1_{t,x} (Q_{\rho_0/4}),$$
where $S_{\ep}, T_{\ep}, R_{\ep}$ are the same as in Lemma \ref{prop_NRST}.
\end{Lemme}
Let us postpone again the proof of Lemma \ref{lemma_Tep} to the Appendix.\\

Now, with these two lemmas at hand we can continue the proof of Proposition \ref{lemme_limit}. Indeed, we start noting that the quantites $\mu_\ep$ and $\eta_\ep$ given in (\ref{Def_MuEspNuEsp}) can be rewritten  with the help of the quantity $N_\ep (\vX, \vY, \vZ)$ defined in Lemma \ref{prop_NRST}. To be more precise, we have 
\begin{equation*}
\begin{aligned}
\mu_\ep =  N_{\ep} (\vu, \vb, \vb) +  N_{\ep} (\vb, \vu, \vu) \quad \mbox{and}\quad
\eta_\ep = N_{\ep} (\vv, \vh, \vh) +  N_{\ep} (\vh, \vv, \vv).
\end{aligned}
\end{equation*}   
Thus, proving the limit \eqref{limit_muetq} turns out to prove that
\begin{equation}\label{limit_ofNep}
\begin{aligned}
& \lim_{\ep \to 0}  N_{\ep} (\vu - \vv, \vb, \vb) + N_{\ep} (\vv, \vb-\vh, \vb) - N_{\ep} (\vv, \vh, \vb-\vh) \\ 
&\hspace{0.5cm}+ N_{\ep} (\vb- \vh, \vu, \vu) + N_{\ep} (\vh, \vu-\vv, \vu) - N_{\ep} (\vh, \vv, \vu-\vv) =0.
\end{aligned}
\end{equation} 
Since $\vu, \vv, \vb, \vh  \in L^\infty_t L^2_x \cap L_{t}^2H_{x}^1 (\Qr)$, we obtain easily that $\vu, \vv, \vb, \vh \in L^3_{t,x}(\Omega)$ by interpolation. By Lemma \ref{prop_NRST}, we know that
$
\lim\limits_{\ep \to 0} N_{\ep}  (\vX, \vY, \vZ) = \lim\limits_{\ep \to 0} (S_{\ep} + T_{\ep} -R_{\ep} ) (\vX, \vY, \vZ)
$ for any $\vX, \vY, \vZ  \in L^3_{t,x}(\Qr)$, so we may replace all the $N_\ep$ by $(S_{\ep} + T_{\ep} -R_{\ep} )$ in the formula above. We have
moreover that $\vu-\vv, \vb - \vh \in L^\infty Lip_x (\Qro)$ from Proposition \ref{Lemma_correctors}, so applying Lemma \ref{lemma_Tep} to each term on the left-hand side of \eqref{limit_ofNep} we can conclude that
$$\lim_{\ep \to 0} \mu_\ep - \eta_\ep = 0,$$
and this completes the proof of Lemma \ref{lemme_limit}\hfill$\blacksquare$\\
\section{Gain of regularity for new variables $\vv$ and $\vh$ - Step 3}\label{subsection_vh}
In this section, we will obtain a regularity result for the new variables $\vv$ and $\vh$ (which satisfy the system (\ref{EquationMHD_companion})) by using the partial regularity theory obtained in \cite{ChCHJ2} and stated in Theorem \ref{Teo_CKNMHD} above. To achieve this task, we only need to check the hypotheses \emph{1)-3)} of this theorem. Indeed, from Proposition \ref{Proposition_vh} and Proposition \ref{Proposition_newPressureForce}, we have already shown that $\vv, \vh \in L^{\infty}_t L^{2}_x \cap L^{2}_t H^{1}_x(\Qro)$, $q, r \in L^{\frac{3}{2}}_{t,x} (\Qro) $, $div(\vk) =0$, $div(\vl) =0$ and $\vk -\vk_0 \in L^{2}_{t,x} (\Qro)$ and $\vl - \vl_0 \in L^{2}_{t,x}(\Qro)$. It still remains to show the suitability of $(\vv, \vh)$, the local information of $\vk$ and $\vl$ and the small gradient condition \eqref{HypothesePetitesseGrad} on $\vv$ and $\vh$.

\begin{itemize}
\item \emph{Suitability of $(\vv, \vh)$}. Let us remark that, from the identity \eqref{regular_ub} and by the dissipativity assumption over $(\vu, \vb)$, we get that the quantity
\begin{eqnarray*}
\lambda&=&-\partial_t(|\vu|^2 + |\vb|^2 )+ \Delta (|\vu|^2 + |\vb|^2 )  - 2 (|\vn \otimes \vu|^2 + |\vn \otimes \vb|^2 )- div ( |\vu|^2 \vb + |\vb|^2 \vu )\\
& &	- 2 \langle div \big(P (\vu + \vb)\big)\rangle  + 2 (\vf \cdot \vu + \vg \cdot\vb) = \lim\limits_{\ep \to 0} \mu_\ep
\end{eqnarray*}
is well-defined as a distribution and moreover we have $\lambda \geq 0$ on $\Qr$. As in addition we have (by Lemma \ref{lemme_limit}) that $\lim\limits_{\ep \to 0} \mu_\ep = \lim\limits_{\ep \to 0} \eta_\ep$, it is easy to find that the quantity
\begin{equation*}
\begin{aligned}
\lim\limits_{\ep \to 0} \eta_\ep
&= -\partial_t(|\vv|^2 + |\vh|^2 )
+ \Delta (|\vv|^2 + |\vh|^2 )
- 2 (|\vn \otimes \vv|^2 + |\vn \otimes \vh|^2 ) - div (|\vv|^2 \vh + |\vh|^2 \vv) \\ 
& \hspace{0.5cm} -2  div ( q \vv + r \vh )
+ 2 (\vv \cdot \vk+ \vh \cdot \vl),
\end{aligned}
\end{equation*}
is a non-negative locally finite Borel measure on $\Qro$, which implies immediately the suitability of $(\vv, \vh)$. 
\item \emph{Local information of $\vk$ and $\vl$}. First, we recall the definition of $ \vk_0$ and $ \vl_0$:
\begin{equation*}
\begin{split}
\vk_0 = - \frac{1}{\Delta} \vn \wedge (\psi \vn \wedge \vf ) \quad \text{and} \quad
\vl_0 = - \frac{1}{\Delta} \vn \wedge (\psi \vn \wedge \vg ),
\end{split}
\end{equation*}
and using the vector calculus identities $ \psi \vn \wedge \vf=\vn \wedge (\psi \vf) - \vn \psi \wedge \vf$ and $\vn \wedge \vn \wedge (\psi \vf)=\vn(\vf \cdot \vn \psi) - \Delta (\psi \vf)$, 
we get that 
\begin{eqnarray*}
\vk_0=- \frac{1}{\Delta} \vn \wedge \left(\vn \wedge (\psi \vf) - \vn \psi \wedge \vf \right) 	=  -  \frac{1}{\Delta} \left( \vn (\vf \cdot \vn \psi)\right) + \psi \vf + \frac{1}{\Delta} \vn \wedge (\vn \psi \wedge \vf) =\vf + \vbe_{\vf},
\end{eqnarray*}
where $\vbe_{\vf}:=  -  \frac{1}{\Delta} \left( \vn (\vf \cdot \vn \psi)\right) + \frac{1}{\Delta} \vn \wedge (\vn \psi \wedge \vf)$. Doing the same computation for $\vl_0$, we get the decomposition
\begin{equation*}
\begin{split}
\vl_0 = \vg + \vga_{\vg} \quad \text{with} \quad \vga_{\vg} := -  \frac{1}{\Delta} \left( \vn (\vg \cdot \vn \psi)\right) + \frac{1}{\Delta} \vn \wedge (\vn \psi \wedge \vg).
\end{split}
\end{equation*}
By the hypotheses \eqref{Formula_MesureDissipativeMHD1},  we have $\vf, \vg \in L^2_{t}H^1_x (\Qr)$ for the forces $\vf, \vg$ in the MHD system \eqref{EquationMHD}. From the proof and the result of Proposition \ref{Lemma_correctors}, we can deduce that $\vbe_{\vf}, \vga_{\vf}$ belong to $ L^2_t Lip_x (\Qro)$. Indeed, using the similar computation as in \eqref{betau_Lip}, we have
\begin{equation*}
\begin{split}
&\|\vbe_{\vf}\|_{L^2_t H^1_x(\Qro)} 
\leq C_{\rho_0} \|\vbe_{\vf}\|_{L^2_t \text{Lip}_x(\Qro)}
\leq C_{\rho_0} \|\vf\|_{L^2_{t,x}(\Qr)} 
\leq C_{\rho_0} \|\vf\|_{L^2_t H^1_x(\Qr)}\\
&\|\vga_{\vg}\|_{L^2_t H^1_x(\Qro)} 
\leq C_{\rho_0} \|\vga_{\vg}\|_{L^2_t \text{Lip}_x(\Qro)}
\leq C_{\rho_0} \|\vg\|_{L^2_{t,x}(\Qr)} 
\leq C_{\rho_0} \|\vg\|_{L^2_t H^1_x(\Qr)}.
\end{split}
\end{equation*}
Thus, $\vk_0, \vl_0 \in L^2_t H^1_x(\Qro) \subset L^2_{t,x} (\Qro)$ and it results in 
$$\vk, \vl \in L^{2}_{t,x} (\Qro) \quad \text{and} \quad \mathds{1}_{\Qro} \vk, \mathds{1}_{\Qro}\vl \in L^2_{t,x} \subset \mathcal{M}_{t,x}^{\frac{10}{7}, 2}.$$
\item \emph{Small gradient condition on $\vv$ and $\vh$}. As we know from Proposition \ref{Lemma_correctors} that $\vbe = \vu -\vv$ and $\vga= \vb -\vh$ belong to $L^\infty_t Lip_x (\Qro)$, thus for any $0<r< \rho_0$, we have 
\begin{equation*}
\begin{split}
\frac{1}{r} \iint_{Q_r}|\vn \otimes \vbe \,(s, y)|^{2} + |\vn \otimes \vga \,(s, y)|^{2} dyds \leq C r^4 \left(\|\vbe\|^2_{L^\infty_t Lip_x (\Qro)} + \|\vga\|^2_{L^\infty_t Lip_x (\Qro)} \right)
\end{split}
\end{equation*}
which results in
\begin{equation*}
\begin{split}
\limsup _{r \rightarrow 0}
\frac{1}{r} \iint_{{Q_r}}|\vn \otimes \vu \,(s, y)|^{2} + |\vn \otimes \vb \,(s, y)|^{2} dyds
= \limsup _{r \rightarrow 0} \frac{1}{r} \iint_{{Q_r}}|\vn \otimes \vv \,(s, y)|^{2} + |\vn \otimes \vh \,(s, y)|^{2} dyds,
\end{split}
\end{equation*}
so the small gradient condition \eqref{HypothesePetitesseGrad} is essentially equal to the following condition: there exists a positive constant $\epsilon^{*}$  which depends only on $\tau_{0} >5$ such that for some $\left(t_{0}, x_{0}\right) \in \Qro$, we have
\begin{equation}\label{smallgradub}
\limsup _{r \rightarrow 0} \frac{1}{r} \iint_{\left]t_{0}-r^{2}, t_{0}+r^{2}\right[ \times B(x_0, r)}|\vn \otimes \vu \,(s, y)|^{2} + |\vn \otimes \vb \,(s, y)|^{2} dyds<\epsilon^{*}. 
\end{equation} 
\end{itemize}
We have now verified all the assumption of Theorem \ref{Teo_CKNMHD} and we can thus conclude a gain of information for the variables $\vv$ and $\vh$:  the couple $(\vv, \vh)$ is H\"older regular (in the time variable and in the space variable) in a neighborhood of $\left(t_{0}, x_{0}\right)$ and from this fact we can easily deduce that 
\begin{equation}\label{gain_reg_vh}
\mathds{1}_{Q_{R_0}} \vv \in \mathcal{M}_{t,x}^{3, \tau_{0}}\qquad \mbox{and}\qquad \mathds{1}_{Q_{R_0}} \vh \in \mathcal{M}_{t,x}^{3, \tau_{0}}, \qquad \tau_0 > 5.
\end{equation}
\section{Gain of regularity for the variables $\vu$ and $\vb$ - Step 4}\label{Sec_GainReguubb}
This section is devoted to the proof of the following statement in which we have a gain of regularity in the original variables $\vu$ and $\vb$. 
\begin{Theoreme}\label{Theorem_main}
Let $\Omega$ be a bounded domain of $\mathbb{R} \times \R$ and $(\vu, P, \vb)$ be a weak solution on $\Omega$ of the MHD equations \eqref{EquationMHD}.
Assume that
\begin{itemize}
\item[1)] $(\vu, \vb, P, \vf, \vg)$ satisfies the conditions :
$$\vu, \vb \in L^{\infty}_t L^{2}_x \cap L^{2}_t \dot{H}^{1}_x(\Omega), \quad \vf, \vg \in L^{2}_{t} H^1_x(\Omega), \quad P \in \mathcal{D}'(\Omega);$$
\item[2)] $(\vu, P, \vb)$ is dissipative in the sense of Definition \ref{Def_dissi};
\item[3)] there exists a positive constant $\epsilon^{*}$  which depends only on $\tau_{0} >5$ such that for some  $\left(t_{0}, x_{0}\right) \in \Omega$, we have
\begin{equation}\label{smallgradub}
\limsup _{r \rightarrow 0} \frac{1}{r} \iint_{\left]t_{0}-r^{2}, t_{0}+r^{2}\right[ \times B(x_0, r)}|\vn \otimes \vu \,(s, y)|^{2} + |\vn \otimes \vb \,(s, y)|^{2} dyds<\epsilon^{*}, 
\end{equation}
\end{itemize}
then there exists a (parabolic) neighborhood $Q_{R_0}$ of $(t_{0}, x_{0})$ such that 
\begin{eqnarray}\label{gain_reg_ub}
\mathds{1}_{Q_{R_0}} \vec{u} \in \mathcal{M}_{t,x}^{3, \tau_{0}}, \quad \mathds{1}_{Q_{R_0}} \vb \in \mathcal{M}_{t,x}^{3, \tau_{0}}.
\end{eqnarray}
\end{Theoreme}
\textbf{Proof.}
Based on the arguments in the Section \ref{subsection_vh} and the conclusion of Proposition \ref{Proposition_vh}, it remains to show that the local information \eqref{gain_reg_vh} on $(\vv, \vh)$ implies the local information \eqref{gain_reg_ub} on $(\vu, \vb)$. 
Since $\tau_0 >3$, by the (generic) local property of Morrey spaces $\|\mathds{1}_{Q_{R_0}}\vf\|_{\mathcal{M}_{t,x}^{3, \tau_0}} \leq C\|\mathds{1}_{Q_{R_0}}\vf\|_{\mathcal{M}_{t,x}^{\tau_0, \tau_0}} =  C\|\vf\|_{L_{t,x}^{ \tau_0}(Q_{R_0})}$ and by the identity $\vu=\vv+\vbe$ (recall (\ref{Def_DifferenceHarmonique})) we obtain that
\begin{equation*}
\begin{split}
\|\mathds{1}_{Q_{R_0}}\vu\|_{\mathcal{M}_{t,x}^{3, \tau_0}} 
\leq \|\mathds{1}_{Q_{R_0}}\vv\|_{\mathcal{M}_{t,x}^{3, \tau_0}}+ \|\mathds{1}_{Q_{R_0}}\vbe\|_{\mathcal{M}_{t,x}^{3, \tau_0}}
\leq \|\mathds{1}_{Q_{R_0}}\vv\|_{\mathcal{M}_{t,x}^{3, \tau_0}}+ \|\vbe\|_{L_{t,x}^{\infty}(Q_{R_0})} < +\infty,
\end{split}
\end{equation*}
where we used the embedding $L_{t,x}^{\infty}(Q_{R_0}) \subset L_{t,x}^{\tau_0}(Q_{R_0})$.
By exactly the same computations, we thus have $\mathds{1}_{Q_{R_0}} \vb \in \mathcal{M}_{t,x}^{3, \tau_{0}}$ and this ends the proof of Theorem \ref{Theorem_main}.
\hfill$\blacksquare$
\section{Regularity of the variables $\vU$ and $\vB$ - Step 5}\label{Sec_RegUB}
This section is devoted to end the proof of Theorem \ref{Theorem_main_original} by making use of Theorem \ref{Theorem_main} and the local regularity result obtained in \cite{ChCHJ1} (see Theorem \ref{Teo_SerrinMHD}).\\

\noindent\textbf{Proof of Theorem \ref{Theorem_main_original}.}
First, since we have the relationships $\vu= \vU + \vB$, $\vb= \vU-\vB$, $\vf= \vF + \vG$ and $\vg = \vF - \vG$, one can check that hypothesis \emph{2)} in Theorem \ref{Theorem_main_original} implies that the quantity $\lambda$
\begin{eqnarray*}
\lambda&=&-\partial_t(|\vu|^2 + |\vb|^2 )+ \Delta (|\vu|^2 + |\vb|^2 )  - 2 (|\vn \otimes \vu|^2 + |\vn \otimes \vb|^2 )- div ( |\vu|^2 \vb + |\vb|^2 \vu )\\
& &	- 2 \langle div \big(P (\vu + \vb)\big)\rangle  + 2 (\vf \cdot \vu + \vg \cdot\vb)
\end{eqnarray*}
is well-defined as a distribution and is a locally finite non-negative measure on $\Omega$, so that  $(\vu, \vb)$ is dissipative solution of equations (\ref{EquationMHD}), i.e., point $2)$ in Theorem \ref{Theorem_main}. Moreover, 
the hypothesis \emph{3)} in Theorem \ref{Theorem_main_original} implies the point $3)$ in Theorem \ref{Theorem_main}. Indeed,
\begin{equation*}
\begin{split}
&\limsup _{r \rightarrow 0} \frac{1}{r}
\iint_{\left]t_{0}-r^{2}, t_{0}+r^{2}\right[ \times B(x_0, r)}|\vn \otimes \vu \,(s, y)|^{2} + |\vn \otimes \vb \,(s, y)|^{2} dyds\\
= & \limsup _{r \rightarrow 0} \frac{1}{r} \iint_{\left]t_{0}-r^{2}, t_{0}+r^{2}\right[ \times B(x_0, r)}|\vn \otimes (\vU + \vB) \,(s, y)|^{2} + |\vn \otimes (\vU - \vB) \,(s, y)|^{2} dyds < \epsilon^*.
\end{split}
\end{equation*}
Applying Theorem \ref{Theorem_main}, we get local information \eqref{gain_reg_ub} on $(\vu, \vb)$, which is indeed equal to
\begin{eqnarray}\label{hy_localubUB}
 \mathds{1}_{Q_{R_0}}\vU, \mathds{1}_{Q_{R_0}}\vB\in \mathcal{M}_{t,x}^{3,\tau_{0}}(\mathbb{R}\times\R)\quad \text{with} \quad \tau_0 >5.
\end{eqnarray} 
Observe now that $(\vU, \vB)$ satisfies the local hypotheses \eqref{LocalHypo1} in Theorem \ref{Teo_SerrinMHD}. We conclude the proof by applying the Theorem \ref{Teo_SerrinMHD} and the arguments stated in the formulas \eqref{ConclusionSerrinMHD}-\eqref{ConclusionSerrinMHD1}.
\hfill$\blacksquare$

\section*{Appendix}

{\bf Proof of Lemma \ref{prop_NRST}.}
For simplicity of the argument, we denote by $\vX_{\ep} := \vX * \theta_\ep$, $\vY_{\ep}: = \vY * \theta_\ep$ and $\vZ_{\ep} := \vZ * \theta_\ep$.
We first remark that, as $div(\vX) = 0$ and $\theta$ is a smooth function on $\R$ with compact support, then for any point $(t,x)$ on the cylinder $Q_{\rho_0/4}$, we have the identity $\displaystyle{\int_{\R}} \vn \theta_{\ep}(y) \cdot\left( \vX (t, x-y) - \vX (t,x) \right) dy = 0$, so that 
$$\int_{\R} \left(\vn \theta_{\ep}(y)\cdot \delta_{y}[\vX ](t, x) \right) (\vY \cdot \vZ )(t, x) d y = 0,$$
which allows us to rewrite $ R_{\ep} (\vX, \vY, \vZ)$ in the following ways : 
\begin{equation}\label{formular_R}
\begin{aligned}
R_{\ep} (\vX, \vY, \vZ)(t,x)
& = \int_{\R} \left(\vn \theta_{\ep}(y)\cdot \delta_{y}[\vX ](t, x)\right) 
(\vY \cdot \vZ )(x-y)d y \\
&\hspace{0.5cm} - \int_{\R} \left(\vn \theta_{\ep}(y)\cdot \delta_{y}[\vX ](t, x)\right)
\left(\vY(x-y) \cdot \vY (x) + \vZ(x-y) \cdot \vZ (x)\right) d y \\
& = div \left(\theta_\ep * [(\vY \cdot \vZ) \vX] \right)
- div \left([\theta_\ep * (\vY \cdot \vZ)] \vX \right)\\
&\hspace{0.5cm} - \int_{\R} \left(\vn \theta_{\ep}(y)\cdot \delta_{y}[\vX ](t, x)\right)
\left(\vY(x-y) \cdot \vY (x) + \vZ(x-y) \cdot \vZ (x)\right) d y, 
\end{aligned}
\end{equation}
where we used the identities 
\begin{equation*}
\begin{aligned}
&  \int_{\R} \left(\vn \theta_{\ep}(y)\cdot \vX ( x-y)\right)
(\vY \cdot \vZ )(x-y)d y =  div \left(\theta_\ep * [(\vY \cdot \vZ) \vX]\right), \\
&  \int_{\R} \left(\vn \theta_{\ep}(y)\cdot \vX ( x)\right)
(\vY \cdot \vZ )(x-y)d y = div \left([\theta_\ep * (\vY \cdot \vZ)] \vX \right).
\end{aligned}
\end{equation*}
In the same manner, for $S_\ep$ and $T_\ep$, we have 
\begin{equation}\label{formular_S}
\begin{aligned}
S_{\ep}(\vX, \vY, \vZ) (t,x)
&=\int_{\R}  \left(\vn \theta_{\ep}(y)\cdot \delta_{y}[\vX ](t, x)\right) \left(\delta_{y}[\vY ](t, x) \right) \cdot\left(\vZ_{\ep}(x)-\vY(x)\right) d y\\
&= \int_{\R}  \left(\vn \theta_{\ep}(y)\cdot \delta_{y}[\vX ](t, x)\right) \left(\vY(x-y)\right) \cdot\left(\vZ_{\ep}(x)-\vY(x)\right) d y\\
& =  \int_{\R}  \left(\vn \theta_{\ep}(y)\cdot \vX ( x-y)\right)
\left(\vY(x-y)\cdot\vZ_{\ep}(x)\right) d y  \\
& \hspace{0.5cm}
- \int_{\R}  \left(\vn \theta_{\ep}(y)\cdot \vX (x)\right)
\left(\vY(x-y) \cdot\vZ_{\ep}(x)\right) d y\\
& \hspace{0.5cm}
- \int_{\R}  \left(\vn \theta_{\ep}(y)\cdot \delta_{y}[\vX ](t, x)\right) \left(\vY(x-y)\cdot\vY(x)\right) d y\\ 
& = \vZ_\ep \cdot \left([(\vX \cdot \vn)\vY ]* \theta_\ep\right) -   \vZ_\ep \cdot \left((\vX \cdot \vn)\vY_\ep\right) \\
& \hspace{0.5cm}
-  \int_{\R}  \left(\vn \theta_{\ep}(y)\cdot \delta_{y}[\vX ](t, x)\right) \left(\vY(x-y)\cdot\vY(x)\right) d y,
\end{aligned}
\end{equation}
and 
\begin{equation}\label{formular_T}
\begin{aligned}
T_{\ep}
& (\vX, \vY, \vZ) (t,x)
= \int_{\R}  \left(\vn \theta_{\ep}(y)\cdot \delta_{y}[\vX ](t, x)\right) \left(\delta_{y}[\vZ ](t, x) \right) \cdot\left(\vY_{\ep}(x)-\vZ(x)\right) d y\\
& = \vY_\ep \cdot \left([(\vX \cdot \vn)\vZ ]* \theta_\ep\right)- \vY_\ep \cdot \left((\vX \cdot \vn)\vZ_\ep\right) 
- \int_{\R}  \left(\vn \theta_{\ep}(y)\cdot \delta_{y}[\vX ](t, x)\right) \left(\vZ(x-y)\cdot\vZ(x)\right) d y.
\end{aligned}
\end{equation}
Gathering the expressions \eqref{formular_R}, \eqref{formular_S} and \eqref{formular_T} together, we have 
\begin{equation*}
\begin{aligned}
(R_{\ep}  - S_{\ep} - T_{\ep})  (\vX, \vY, \vZ)
& =  div \left(\theta_\ep * [(\vY \cdot \vZ) \vX] \right)
- div \left([\theta_\ep * (\vY \cdot \vZ)] \vX \right)- \vZ_\ep \cdot \left([(\vX \cdot \vn)\vY ]* \theta_\ep\right)
\\
& \hspace{0.5cm} -  \vY_\ep \cdot \left([(\vX \cdot \vn)\vZ ]* \theta_\ep\right)  + \vZ_\ep \cdot \left((\vX \cdot \vn)\vY_\ep\right)
+ \vY_\ep \cdot \left((\vX \cdot \vn)\vZ_\ep\right)\\
& =  div \left(\theta_\ep * [(\vY \cdot \vZ) \vX] 
- [\theta_\ep * (\vY \cdot \vZ)] \vX \right)\\
& \hspace{0.5cm}- \vZ_\ep \cdot \left([(\vX \cdot \vn)\vY ]* \theta_\ep\right) -  \vY_\ep \cdot \left([(\vX \cdot \vn)\vZ ]* \theta_\ep\right)  +div \left((\vY_\ep \cdot \vZ_\ep) \vX\right),
\end{aligned}
\end{equation*}
where in the last line we used the identity 
$\vZ_\ep \cdot \left((\vX \cdot \vn)\vY_\ep\right) + \vY_\ep \cdot \left((\vX \cdot \vn)\vZ_\ep\right) = div \left((\vY_\ep \cdot \vZ_\ep) \vX\right)$ due to divergence free assumption on $\vX$. Recall now the expression of $N_{\ep} (\vX, \vY, \vZ)$ given in (\ref{four_quantities}): 
$$N_{\ep} (\vX, \vY, \vZ)=\vY_{\ep} \cdot\left([(\vX \cdot \vn ) \vZ] * \theta_{\ep}\right) + \vZ_{\ep} \cdot\left([(\vX \cdot \vn ) \vY] * \theta_{\ep}\right)-div((\vY \cdot \vZ) \vX),$$
we thus have
\begin{equation*}
\begin{aligned}
(R_{\ep}  - S_{\ep} - T_{\ep} + N_\ep) (\vX, \vY, \vZ)
& =  div \left(\theta_\ep * [(\vY \cdot \vZ) \vX] 
- [\theta_\ep * (\vY \cdot \vZ)] \vX \right) + div \left((\vY_\ep \cdot \vZ_\ep) \vX -(\vY \cdot \vZ) \vX \right).
\end{aligned}
\end{equation*}
Since $\vX, \vY, \vZ \in L^{3}_{t,x} (Q_{\rho_0})$, we have $\lim\limits_{\ep \to 0} \left( \theta_\ep * [(\vY \cdot \vZ) \vX] - [\theta_\ep * (\vY \cdot \vZ)] \vX \right) = 0$  in $L^1_{t,x} (Q_{\rho_0/4} )$, hence, 
$$ \lim\limits_{\ep \to 0} div \left(\theta_\ep * [(\vY \cdot \vZ) \vX] 
- [\theta_\ep * (\vY \cdot \vZ)] \vX \right) = 0 \quad \text{in} \quad \mathcal{D}'(Q_{\rho_0/4}).$$
Moreover, we have limit $\lim\limits_{\ep \to 0} div \left((\vY_\ep \cdot \vZ_\ep) \vX -(\vY \cdot \vZ) \vX \right) = 0$ in $\mathcal{D}'( Q_{\rho_0/4} )$,
which allows us to conclude that 
$$\lim_{\ep \to 0} (N_{\ep} + R_{\ep}- S_{\ep} - T_{\ep}) (\vX, \vY, \vZ)= 0.$$
The lemma \ref{prop_NRST} is now proved.\hfill$\blacksquare$\\

\noindent{\bf Proof of Lemma \ref{lemma_Tep}.}
By the definitions of $S_{\ep},T_{\ep}, R_{\ep}$ given in (\ref{four_quantities}), we can rewrite the operator $(S_{\ep} + T_{\ep} - R_{\ep}) (\vX, \vY, \vZ)$ in the following manner:
\begin{equation*}
\begin{aligned}
&(S_{\ep} + T_{\ep} - R_{\ep}) (\vX, \vY, \vZ)(t,x) \\
& =  \int_{\R} \left(\vn \theta_{\ep}(y)\cdot \delta_{y}[\vX ](t, x)\right) \delta_{y}[\vY ](t, x)  d y \cdot \left((\vZ * \theta_\ep)(t, x) - \vZ(t, x) + \vZ(t, x)-\vY(t, x)\right) \\
&\hspace{0.5cm} + \int_{\R} \left(\vn \theta_{\ep}(y)\cdot \delta_{y}[\vX ](t, x)\right)\delta_{y}[\vZ ](t, x)   d y \cdot \left( (\vY * \theta_\ep)(t, x) - \vY(t, x) + \vY(t, x)-\vZ(t, x)\right)\\
&\hspace{0.5cm} -  \int_{\R} \left(\vn \theta_{\ep}(y)\cdot \delta_{y}[\vX ](t, x)\right)\left(\delta_{y}[\vY ](t, x) +  \vY(t, x)-\vZ(t, x)\right) \cdot \left(\delta_{y}[\vZ](t, x) + \vZ (t, x) - \vY (t, x)\right)d y\\
&=  \int_{\R} \left(\vn \theta_{\ep}(y)\cdot \delta_{y}[\vX ](t, x)\right) \delta_{y}[\vY ](t, x)  d y \cdot \left((\vZ * \theta_\ep)(t, x) - \vZ(t, x)\right) \\
&\hspace{0.5cm} +  \int_{\R} \left(\vn \theta_{\ep}(y)\cdot \delta_{y}[\vX ](t, x)\right)\delta_{y}[\vZ ](t, x)   d y \cdot \left( (\vY * \theta_\ep)(t, x) - \vY(t, x)\right)\\
&\hspace{0.5cm} - \int_{\R} \left(\vn \theta_{\ep}(y)\cdot \delta_{y}[\vX ](t, x)\right)\left(\delta_{y}[\vY ](t, x)\right) \cdot \left(\delta_{y}[\vZ](t, x)\right)dy,
\end{aligned}
\end{equation*}
where in the last line we used the identity $\displaystyle{\int_{\R}} \vn \theta_{\ep}(y)\cdot \delta_{y}[\vX ](t, x) d y = 0$ due to divergence free condition of $\vX$. Next, recall the definition of the mollifier
$\theta_\ep = \frac{1}{\ep^3} \theta (\frac{x}{\ep})$, we have $\displaystyle{\int_{\R}}\theta_\ep dx= 1$, $\text{supp} (\theta_\ep) \subset B(0, \ep)$ and $|\vn \theta_{\ep}(y)| \leq \frac{C}{\ep^4}$, hence,
\begin{equation}\label{bigformula_Tep}
\begin{aligned}
\vert(S_{\ep} + T_{\ep} - R_{\ep}) (\vX, \vY, \vZ)(t,x) \vert
&\leq C \frac{1}{\ep^4} \int_{|y| < \ep} \left|\delta_{y}[\vX](t, x)\right|\,\left|\delta_{y}[\vY ](t, x)\right|\, d y \left|(\vZ * \theta_\ep)(t, x)-\vZ(t, x)\right|\\
&\hspace{0.5cm} + C \frac{1}{\ep^4} \int_{|y| < \ep} \left|\delta_{y}[\vX](t, x)\right|\, \left|\delta_{y}[\vZ ](t, x)\right|  d y \, \left|(\vY * \theta_\ep)(t, x)-\vY(t, x)\right|\\
&\hspace{0.5cm} + C \frac{1}{\ep^4} \int_{|y| < \ep} \left| \delta_{y}[\vX](t, x)\right| \,\left|\delta_{y}[\vY ](t, x)\right| \, \left|\delta_{y}[\vZ](t, x)\right|  d y\\
& \leq  C \frac{1}{\ep^7} \int_{|y| < \ep} \left|\delta_{y}[\vX](t, x)\right|\,\left|\delta_{y}[\vY ](t, x)\right|\, d y \times \int_{|z| < \ep} \left| \delta_{z}[\vZ ](t, x)\right| d z\\
&\hspace{0.5cm} + C \frac{1}{\ep^7} \int_{|y| < \ep} \left|\delta_{y}[\vX](t, x)\right|\, \left|\delta_{y}[\vZ ](t, x)\right|  d y \times \int_{|z| < \ep} \left| \delta_{z}[\vY](t, x)\right| d z\\
&\hspace{0.5cm} + C \frac{1}{\ep^4} \int_{|y| < \ep} \left| \delta_{y}[\vX](t, x)\right| \,\left|\delta_{y}[\vY ](t, x)\right| \, \left|\delta_{y}[\vZ](t, x)\right|  d y
\end{aligned}
\end{equation}
where we used the following estimates 
$$ \left|(\vZ * \theta_\ep)(t, x)-\vZ(t, x)\right| = \left| \int_{\R} \theta_\ep (z) \left(\vZ(t, x-z) -\vZ(t, x)  \right) d z  \right|  \leq C \frac{1}{\ep^3} \int_{|z| < \ep} \left| \delta_{z}[\vZ ](t, x)\right| d z,$$
$$\left|(\vY * \theta_\ep)(t, x)-\vY(t, x)\right| = \left| \int_{\R} \theta_\ep (z) \left(\vY(t, x-z) -\vY(t, x)\right) d z  \right| \leq C \frac{1}{\ep^3} \int_{|z| < \ep} \left| \delta_{z}[\vY ](t, x)\right| d z.$$
Now, using Hölder's inequality to the right-hand side of  \eqref{bigformula_Tep}, we have  
\begin{equation*}
\begin{aligned}
\vert(S_{\ep} + T_{\ep} - R_{\ep}) (\vX, \vY, \vZ)(t,x) \vert
& \leq C \frac{1}{\ep^4} \|\delta_{\cdot}[\vX](t, x)\|_{L_y^3(B(0,\ep))}  \|\delta_{\cdot}[\vY](t, x)\|_{L_y^3(B(0,\ep))}   \|\delta_{\cdot}[\vZ](t, x)\|_{L_y^3(B(0,\ep))} .
\end{aligned}
\end{equation*}
Let us turn to prove the limit in $L^1_{t,x} (Q_{\rho_0/4})$ norm. Taking the $L^1$ norm of $(S_{\ep} + T_{\ep} - R_{\ep})(\vX, \vY, \vZ)(t,x)$ in time and in space variable over the cylinder $Q_{\rho_0/4}$ and applying Hölder's inequality we obtain 
\begin{eqnarray}
&&\int_{Q_{\rho_0/4}}  \vert (S_{\ep} + T_{\ep} - R_{\ep})(\vX, \vY, \vZ)(t,x) \vert d t d x\notag\\ 
&\leq & \; C \frac{1}{\ep^4}
\int_{Q_{\rho_0/4}} \|\delta_{\cdot}[\vX](t, x)\|_{L_y^3(B(0,\ep))}  \|\delta_{\cdot}[\vY](t, x)\|_{L_y^3(B(0,\ep))}   \|\delta_{\cdot}[\vZ](t, x)\|_{L_y^3(B(0,\ep))} d t d x\notag\\
&\leq & \;C \frac{1}{\ep^4}
\left(\int_{Q_{\rho_0/4}} 
\|\delta_{\cdot}[\vX](t, x)\|^3_{L_y^3(B(0,\ep))} \right)^{\frac{1}{3}}
\left(\int_{Q_{\rho_0/4}} \|\delta_{\cdot}[\vY](t, x)\|^3_{L_y^3(B(0,\ep))} \right)^{\frac{1}{3}}\notag\\
&& \times \left(\int_{Q_{\rho_0/4}} \|\delta_{\cdot}[\vZ](t, x)\|^3_{L_y^3(B(0,\ep))}\right)^{\frac{1}{3}}.\label{integ_Tep}
\end{eqnarray}
Recalling that $\vX, \vY, \vZ \in L_{t,x}^3 (\Qro)$ and at least one of them belongs to $L_t^\infty Lip_x (\Qro)$, so we may suppose for instance that $\vX \in L_t^\infty Lip_x (\Qro) \cap L_{t,x}^3 (\Qro)$ and $\vY, \vZ \in L_{t,x}^3 (\Qro)$. For the term involving $\vX$ on the right-hand side of \eqref{integ_Tep}, by a change of variables, we get
\begin{equation*}
\begin{aligned}
\int_{Q_{\rho_0/4}} 
\|\delta_{\cdot}[\vX](t, x)\|^3_{L_y^3(B(0,\ep))} d t d x
= \int_{Q_{\rho_0/4}}  \int_{|y|<\ep} |\vX (t, x-y) - \vX (t,x) |^3 d y d t d x
\leq C \ep^6 \| \vn \otimes \vX \|^3_{L^\infty_{t,x}} |Q_{\rho_0/4}|,
\end{aligned}
\end{equation*}
and thus we have 
\begin{equation}\label{firstX}
\frac{1}{\ep^2}\left(\int_{Q_{\rho_0/4}} \|\delta_{\cdot}[\vX](t, x)\|^3_{L_y^3(B(0,\ep))} d t d x\right)^{\frac{1}{3}}
\leq C \| \vn \otimes \vX \|_{L^\infty_{t,x}} |Q_{\rho_0/4}|^{\frac{1}{3}}.
\end{equation}
Moreover, for $\vY, \vZ \in L_{t,x}^3 (\Qro)$, by a change of variables, we have
\begin{equation}\label{lasttwoYZ}
\frac{1}{\ep} \left(\int_{Q_{\rho_0/4}} \|\delta_{\cdot}[\vY](t, x)\|^3_{L_y^3(B(0,\ep))} \right)^{\frac{1}{3}}  \to 0,  \quad \frac{1}{\ep} \left(\int_{Q_{\rho_0/4}} \|\delta_{\cdot}[\vZ](t, x)\|^3_{L_y^3(B(0,\ep))}\right)^{\frac{1}{3}}  \to 0, \quad \text{as} \quad \ep \to 0.
\end{equation} 
Note that the convergence holds here by the dominated convergence theorem. Substituting the estimates \eqref{firstX} and \eqref{lasttwoYZ} into \eqref{integ_Tep}, we conclude that
$$
\lim_{\ep \to 0}(S_{\ep} + T_{\ep} - R_{\ep})(\vX, \vY, \vZ) = 0 \quad \text{in} \quad L^1_{t,x} (Q_{\rho_0/4}).
$$
The Lemma \ref{lemma_Tep} is now proved.
\hfill$\blacksquare$\\

\noindent{\bf Acknowledgments:} J. \textsc{He} is supported by the \emph{Sophie Germain} Post-doc program of the \emph{Fondation Math\'ematique Jacques Hadamard}.

\end{document}